\documentclass[twoside]{article}
\usepackage{amssymb,amsfonts}
\usepackage{graphicx}

\newtheorem{theorem}{Theorem}[section]        
\newtheorem{definition}[theorem]{Definition}     
\newtheorem{lemma}[theorem]{Lemma}             
\newtheorem{remark}[theorem]{Remark}           

\newtheorem{proposition}[theorem]{Proposition}
\newenvironment{proof}{{\it Proof.}}{\mbox{\ }\hfill 
$\square$\vskip3mm}         
\newcommand{\R}{{\mathbb R}}
\newcommand{\C}{{\mathbb C}}
\renewcommand{\Re}{\,{\rm Re}\,}
\renewcommand{\Im}{\,{\rm Im}\,}
\renewcommand{\epsilon}{\varepsilon}
\newcommand{\N}{{\mathbb N}}

\newcommand{\eqnref}[1]{{\rm (\ref{#1})}}
\newcommand{\mysection}[2]{\section*{{\large\bf #1. #2}}%
\setcounter{equation}{0}\setcounter{section}{#1}%
\setcounter{theorem}{0}}

\newcommand{\mylabel}[1]{\def\nummer{#1}
\label{\nummer}}

\advance\textheight-\headheight      
\advance\textheight-\headsep          
\advance\textheight by \topskip
\skip\footins 9pt plus 4pt minus 2pt
\begin{document}
\thispagestyle{empty}
\centerline{\bf M.~V.~Keldysh Institute of Applied Mathematics} 
\centerline{\bf Russian Academy of Sciences}
\vskip 2.5cm
\centerline{\bf R. Denk, R. Mennicken, and L. R. Volevich}
\vskip 1.7cm
\centerline{\bf BOUNDARY VALUE PROBLEMS FOR A CLASS OF}
\vskip 0.3cm
\centerline{\bf  ELLIPTIC OPERATOR PENCILS }
\vskip 2 cm
\centerline{Preprint 58 (1998)}
\vskip 5cm
\centerline{\bf Moscow 1998}

\newpage
\thispagestyle{empty}
\centerline{Abstract}
R. Denk, R. Mennicken, and L. R. Volevich\footnote{Supported in part 
by the Deutsche Forschungsgemeinschaft and by Russian Foundation of 
Fundamental Research, Grant 97-01-00541}.
Boundary value problems for a class of elliptic operator pencils.

\medskip
In this paper operator pencils $A(x,D,\lambda)$ are studied which act
on a manifold with boundary and satisfy the condition of
$N$-ellipticity with parameter, a generalization of the notion of
ellipticity with parameter as introduced by Agmon and
Agranovich--Vishik. Sobolev spaces corresponding to a Newton polygon
are defined and investigated; in particular it is possible to describe 
their trace spaces. With respect to these spaces, an a priori estimate 
holds for the Dirichlet boundary value problem connected with an
$N$-elliptic pencil, and a right parametrix is constructed.

\bigskip



\newpage

\mysection{1}{Introduction}

In this paper we consider operator pencils of the form
\begin{equation}\mylabel{1-1}
A(x,D,\lambda) = A_{2m}(x,D)+\lambda A_{2m-1}(x,D)
+ \cdots + \lambda^{2m-2\mu} A_{2\mu}(x,D)
\end{equation}
acting on a smooth manifold $M$ with smooth boundary $\Gamma$.
 Here $m$ and $\mu$ are integer numbers with $m>\mu\ge 0$, 
$A_{2\mu},\ldots,
A_{2m}$ are partial differential operators in $M$ with infinitely smooth
coefficients and  $\lambda$ is a complex parameter. We assume that
\begin{equation}\mylabel{1-2}
 A_j(x,D) = \sum_{|\alpha|\le j} a_{\alpha j}(x) D^\alpha \quad (j=2\mu,\,
2\mu+1, \ldots, 2m)
\end{equation}
is a differential operator of order $j$ with scalar coefficients $a_{\alpha 
j}(x)\in C^\infty(\overline M)$. 
As usual, we use for multi-indices $\alpha=(\alpha_1,\ldots,
\alpha_n)$ the notation
\begin{equation}\mylabel{1-3}
  D^\alpha = D_1^{\alpha_1}\cdots D_n^{\alpha_n}\,,\quad 
D_j = -i \frac{\partial}{\partial x_j}\,,\quad |\alpha| =
\alpha_1 + \ldots + \alpha_n\,.
\end{equation}

There is a close connection between pencils of the form \eqnref{1-1}
and general parabolic problems which we want to describe briefly. An
important tool in the field of elliptic and parabolic problems is the
concept of the Newton polygon. For a given polynomial
\begin{equation}\mylabel{1-3aa}
 P(\xi,\lambda) = \sum_{\alpha,k} a_{\alpha k} \xi^\alpha
\lambda^k\,,
\end{equation}
where $\xi=(\xi_1,\ldots,\xi_n)\in\R^n$ and $\lambda\in\C$, let
$\nu(P)$ be the set of all integer points $(i,k)$ such that an
$\alpha$ exists with $|\alpha|=i$ and $p_{\alpha k}\not=0$. Then the
Newton polygon $N(P)$ is defined as the convex hull of all points in
$\nu(P)$, their projections on the coordinate axes and the origin. The 
polynomial $P(\xi,\lambda)$ is called $N$-parabolic (see
\cite{gindikin-volevich}, Chapter 2) if $N(P)$ has no
sides parallel to the coordinate axes and if the inequality
\begin{equation}\mylabel{1-3a}
|P(\xi,\lambda)|>\delta \sum_{(i,k)\in N(P)} |\xi|^i \, |\lambda|^k
\end{equation}
holds for all $\lambda\in\C$ with $\Im \lambda < \lambda_0$ where
$\delta > 0$ and $\lambda_0$ are constants. An important example of
such polynomials is the product of polynomials
$P_1(\xi,\lambda),\ldots, P_N(\xi,\lambda)$ which are
quasi-homogeneous and $2b_j$-parabolic in the sense of Petrovskii with 
different weights $2b_j$ ($j=1,\ldots,N$). Note that in this case
$P(\xi,\lambda)$ is no quasi-homogeneous function in
$(\xi,\lambda)$. Similarly (see \cite{denk-mennicken-volevich}),
 the polynomial $P(\xi,\lambda)$ is called
$N$-elliptic with parameter along some ray ${\cal L}$ of the complex
plane if \eqnref{1-3a} holds for all $\xi\in\R^n$ and all $\lambda\in
{\cal L}\,,\;|\lambda|\ge R\,,$ with large enough $R$. 
This type of polynomials
appears, for instance, if one considers Douglis--Nirenberg systems
$A(x,D)-\lambda I$ which are elliptic with parameter. On manifolds
without boundary Douglis--Nirenberg systems were investigated by
Kozhevnikov \cite{kozhevnikov} and by the authors in
\cite{denk-mennicken-volevich}. It turned out that an equivalent
condition for unique solvability of such a system and sharp a
priori estimate is the condition that for every $x$ the determinant
\[ P(x,\xi,\lambda) = \det (A(x,D)-\lambda I) \]
satisfies inequality \eqnref{1-3a}.

The basic idea of the Newton polygon method for the problems mentioned 
above is to assign to $\lambda$ various weights $r_j$ which are
defined by the Newton polygon. For each of these weights we obtain a
different principal part of $P(\xi,\lambda)$ which we denote by
$P_{r_j}(\xi,\lambda)$. On a manifold without boundary there is a finite
open covering $\{U_j\}_j$ of the set of all $(\xi,\lambda)$ and a
corresponding partition of unity $\sum_j \psi_j(\xi,\lambda)\equiv 1$
such that $P(D,\lambda)\psi_j(D,\lambda)$ differs from the
corresponding principal part $P_{r_j}(D,\lambda)\psi_j(D,\lambda)$
only by a small regular perturbation. This allows estimates and
existence results for the operators $P(D,\lambda)$,
cf. \cite{denk-mennicken-volevich} for $N$-elliptic systems and
\cite{gindikin-volevich} for parabolic problems.

Now let us consider the same problems on a manifold with boundary. For 
instance, let $P(D,\lambda)$ be the product of two operators which are 
parabolic in the sense of Petrovskii, i.e.
\[ P(D,\lambda) = (\lambda+A_{2p}(D))\, (\lambda+A_{2q}(D))\,,\]
where $\lambda+A_{2p}(D)$ and $\lambda+A_{2q}(D)$ are $2p$- and
$2q$-parabolic operators, respectively, with $p>q$. If, for instance,
we assign to $\lambda$ the weight $r_1=2q$, we obtain the principal
part $P_{r_1}(D,\lambda) = A_{2p}(D) A_{2q}(D) + \lambda A_{2p}(D)$
which is of the form \eqnref{1-1}. If we take the weight $r_2$ with
$2q<r_2<2p$ the corresponding principal part is $P_{r_2}(D,\lambda) =
\lambda A_{2p}(D)$. The operator $P_{r_2}(D,\lambda)$ has to be
supplied with $p$ boundary conditions while the operator
$P(D,\lambda)$ needs $p+q$ boundary conditions. Thus we can see here
that $P(D,\lambda)$ is (after dividing by $\lambda$) a singular
perturbation of the principal part $P_{r_2}(D,\lambda)$. A similar
situation occurs if the weight of $\lambda$ is larger than $2p+2q$.

So we see from this example that operator pencils of the form
\eqnref{1-1} and singular perturbations naturally arise when we deal
with $N$-parabolic problems on manifolds with boundary. If we consider 
boundary value problems elliptic in the sense of Douglis--Nirenberg,
 the situation is the same
or even more complicated.

As a first step to handle these problems we consider as a model
problem operator pencils of the form \eqnref{1-1} and the
corresponding Dirichlet problem on manifolds with boundary. The aim of 
this paper is to show that the Newton polygon provides an easy and
direct approach to the Sobolev spaces where the pencil acts and to the 
proof of a priori estimate. In particular, we obtain 
a description of the trace spaces which is formulated in the general
context of Sobolev spaces corresponding to Newton polygons. We hope to 
study in a subsequent paper boundary value problems for
 general $N$-parabolic operators on manifolds with
boundary.

For $r,s\in\R$ let the Sobolev space $H^{(r,s)}(\R^n)$ be defined by
\begin{eqnarray}\mylabel{1-4}
\lefteqn{ H^{(r,s)} (\R^n) =}\\
 & & \{ u\in {\cal S}'(\R^n): (|\xi|^2+1)^{s/2}
(|\xi|^2 + |\lambda|^2)^{(r-s)/2} Fu(\xi) \in L_2(\R^n)\}\nonumber
\end{eqnarray}
where $Fu$ denotes the Fourier transform of $u$. The norm in $H^{(r,s)} 
(\R^n)$ is given by 
\begin{equation}\mylabel{1-4a}
\| u\|_{r,s} := \left( \int_{\R^n}  (|\xi|^2+1)^s
(|\xi|^2 + |\lambda|^2)^{r-s} |Fu(\xi) |^2 \,d\xi\right)^{\frac 1 2}\,.
\end{equation}
Restricting the distributions
belonging to $H^{(r,s)}(\R^n)$ to the right half space $\R^n_+ := \{ x =
(x',x_n)\in\R^n: x_n>0\}$, we obtain the Sobolev space $H^{(r,s)}(\R^n_+)$.
See Section 2 for the description of the norm in this Sobolev space.
In the standard way we can also define $H^{(r,s)}(M)$ using local 
coordinates.

For every $r$ and $s$ 
the operator pencil \eqnref{1-1} acts continuously from  $H^{(r,s)}$ to 
$H^{(r-2m,s-2\mu)}$. In what follows in connection with the Dirichlet
problem for \eqnref{1-1}  we will restrict ourselves to the case $r=m, s=\mu$,
i.e. we will realize \eqnref{1-1} as an operator from $H^{(m,\mu)}$ onto
$H^{(-m,-\mu)}$. We will assume this pencil to be 
elliptic with parameter along the ray $[0,\infty)$
 in the following sense: Denote by
\begin{equation}\mylabel{1-6}
A_j^{(0)}(x,\xi) := \sum_{|\alpha|=j} a_{\alpha j}(x)\xi^\alpha\quad (j=2\mu,
\ldots,2m)
\end{equation}
the principal symbol of $A_j$, where $\xi^\alpha = \xi_1^{\alpha_1}\cdots
\xi_n^{\alpha_n}$ for $\xi=(\xi_1,\ldots,\xi_n) = (\xi',\xi_n)$, and by
\begin{equation}\mylabel{1-7}
A^{(0)}(x,\xi,\lambda) := A_{2m}^{(0)}(x,\xi) + \lambda A_{2m-1}^{(0)}(x,
\xi) + \ldots + \lambda^{2m-2\mu}A_{2\mu}^{(0)}(x,\xi)
\end{equation}
the principal symbol of $A(x,D,\lambda)$. Then our main assumption is
that the estimate
\begin{equation}\mylabel{1-8}
| A^{(0)}(x,\xi,\lambda)| \ge C |\xi|^{2\mu}\,(\lambda + |\xi|)^{2m-2\mu}
\quad (\xi\in\R^n, \lambda\in [0,\infty), x\in \overline M)
\end{equation}
holds where the constant $C$ does not depend on $x,\xi$ or $\lambda$. In the
case $\mu=0$ this is the usual definition of ellipticity with
parameter which was introduced by Agmon \cite{agmon} and
Agranovich--Vishik \cite{agranovich-vishik}. 
 Therefore we may assume in the
following that $\mu>0$. In this case even for $\lambda\not=0$
the principal symbol $A^{(0)}(x,\xi,\lambda)$ vanishes for $\xi=0$
 which causes the main 
difficulties in proving existence results and estimates. Note that
the symbol $A^{(0)}(x,\xi,\lambda)$ is homogeneous in $\xi$ and $\lambda$
of degree $2m$, as it is the case for the problems treated in 
\cite{agranovich-vishik}.

We will consider boundary value problems in $\R^n$, $\R^n_+$ and $M$. For this
we will describe the space of 
traces of functions $u\in H^{(m,\mu)}(\R^n_+)$, i.e. the space
\begin{equation}\mylabel{1-9}
\{ D_n^{j-1} u(x',0): u\in H^{(m,\mu)}(\R^n_+)
\}\quad \mbox{ for } j=1,\ldots,m\,.
\end{equation}
This will be done in a more general context in Section 2 where Sobolev
spaces corresponding to Newton polygons are considered. The space 
$H^{(m,\mu)}(\R^n_+)$ appears to be a special case of the space
$H^\Xi(\R^n_+)$ where $\Xi(\xi,\lambda)$ is the weight function
corresponding to the Newton polygon $N(P)$ of a polynomial $P(\xi,
\lambda)$ in $\xi\in\R^n$ and $\lambda\in\C$. It turns out that
 the trace space $\{ D_n^{j-1} u(x',0):
u\in H^\Xi(\R^n_+)\}$ is given by $H^{\Xi^{(-j+\frac 1 2)}}(\R^{n-1})$
where $\Xi^{(-j+\frac 1 2)}(\xi',\lambda)$ denotes the weight function
corresponding to the Newton polygon which is constructed from $N(P)$
by a shift of length $j-\frac 1 2$ to the left. Cf. Section 2 for details.
In particular, in the case of the operator pencil \eqnref{1-1}
the trace spaces have the form $H^{(m_j,\mu_j)}(\R^{n-1})$ where
the parameters $m_j$ and $\mu_j$ can be seen directly from the 
corresponding Newton polygon.

In Section 5 we consider  the Dirichlet boundary
problem
\begin{eqnarray}
A(x, D,\lambda)\, u(x)  & = & f(x) \quad\,\mbox{in } M\,,
\mylabel{1-10}\\
\Big(\frac{\partial}{\partial\nu}\Big)^{j-1}\,
 u(x) & = & g_j(x)  \quad(j=1,\ldots,m)\;\mbox{ on }\Gamma
\,.\mylabel{1-11}
\end{eqnarray}
where $\frac{\partial}{\partial\nu}$ denotes the derivative in the
direction of the inner normal to the boundary. The main theorem states
that for every solution $u\in H^{(m,\mu)}(M)$ of the boundary value 
problem \eqnref{1-10}--\eqnref{1-11}
 the a priori estimate
\begin{equation}\mylabel{1-12}
\| u \|_{m,\mu} \le C \Big( \|f\|_{-m,-\mu} + \sum_{j=1}^m \|g_j\|_{m_j,
\mu_j}+ \lambda^{m-\mu} \|u\|_{L_2(M)}\Big)
\end{equation}
holds for $\lambda\ge \lambda_0$ with a constant $C$ not depending 
on $\lambda$ or $u$. The proof
of this theorem is essentially based on estimates of the solution of
an ordinary differential equation which arises from 
\eqnref{1-10}--\eqnref{1-11} by fixing $x\in \Gamma$, rewriting the
boundary value problem in coordinates corresponding to $x$ and taking
the partial Fourier transform with respect to the first $n-1$
variables. Estimates for the fundamental solution of the resulting
ordinary differential equation can be found in Section 4 and use the
precise knowledge about the zeros of the principal
 symbol $A^{(0)}(x,\xi,\lambda)$
considered as a polynomial in $\xi_n$.

These zeros can (for large $\lambda$) be arranged in two groups,
one group remaining bounded for $\lambda\to\infty$, the other
group of zeros being exactly of order $O(\lambda)$ for  $\lambda\to\infty$.
To obtain this result we have to impose an additional condition on
the principal symbol $A^{(0)}(x,\xi,\lambda)$ which is the analogue
of the condition of regular degeneration which is known from the
theory of singular perturbations (cf. Vishik-Lyusternik
\cite{vishik-lyusternik}).
The details can be found in Section 3.

As mentioned above, there is a close connection between pencils of the 
form \eqnref{1-1} and elliptic boundary value problems with small
parameter. Nazarov obtained in \cite{nazarov} a priori estimates under 
the assumption that the fundamental solutions fulfill some estimates
which are 
similar to those proved in Section 4 below. (The norms used in
\cite{nazarov} differ slightly from the norms used in the present
paper.) In several papers Frank and other authors investigated
singular perturbed problems and corresponding a priori estimates,
cf. \cite{frank} and the references therein. The use of the Newton
polygon method which gives the connection to general parabolic
problems as described above, seems to be new even for singular
perturbed problems.

\mysection{2}{Newton's polygon and functional spaces\\ corresponding to it}

In this section we consider a polynomial $P(\xi,\lambda)$ of the
form \eqnref{1-3aa} and its Newton polygon $N(P)$ which was defined in
the Introduction.
 For a detailed discussion of the
Newton polygon, we refer the reader to Gindikin-Volevich
\cite{gindikin-volevich},  Chapters 1 and 2.

To construct function spaces corresponding to the Newton polygon, we
consider the weight function
\begin{equation}\mylabel{3-2}
\Xi_P(\xi,\lambda) := \sum_{(i,k)\in N(P)} |\xi|^i\,|\lambda|^k\,,
\end{equation}
where the summation on the right-hand side is extended over all 
integer points of $N(P)$. The Sobolev space $H^\Xi$ will arise as
a special case of the following more general definition which is taken
from Volevich-Paneah \cite{volevich-panejah}. It can be seen directly 
that the function $\sigma(\xi) := \Xi_P(\xi,\lambda)$ satisfies the
condition which appears in this definition (cf. also Remark
\ref{3.0a} below). In the following, the Fourier transform $F$ is defined
by 
\begin{equation}\mylabel{3-2a}
Fu(\xi) = \frac{1}{(2\pi)^{\frac n 2}} \int_{\R^n} e^{-ix\cdot \xi}
u(x) \, d x
\end{equation}
for $u\in {\cal S}(\R^n)$, the definition is extended in the usual way 
to distributions $u\in {\cal S}'(\R^n)$.

\begin{definition}\mylabel{3.0}{\rm  Let $\sigma(\xi)$ be a 
continuous function on $\R^n$ with values in $\R_+$ and assume
that $\sigma(\xi)\sigma^{-1}(\eta)\le C(1+|\xi-\eta|^N)$ holds
for all $\xi,\eta\in\R^n$ with constants $C$ and $N$ not depending
on $\xi$ or $\eta$.
 Then $H^\sigma$ is defined as the space
of all distributions $u$ in ${\cal S}'(\R^n)$ such that $\sigma(\xi)
Fu(\xi)\in L_2(\R^n)$. The space $H^\sigma$ is endowed with the norm
\begin{equation}\mylabel{3-3}
\|u\|_{\sigma,\R^n} := \Big( \int_{\R^n} \sigma^2(\xi) |Fu(\xi)|^2\,d\xi
\Big)^{1/2} \,.
\end{equation}
}
\end{definition}

\begin{proposition}\mylabel{3.1} {\rm (See Volevich-Paneah  
\cite{volevich-panejah}.)} Let $\sigma(\xi,\lambda)$ be a continuous
function of $\xi$ and assume that 
\[ \sigma(\xi,\lambda)\sigma^{-1}(\eta,\lambda) \le C_1 (1+|\xi-\eta|^N)\]
holds with a constant $C_1$ not depeding on $\xi,\,\eta$ or
$\lambda$. Let
\[ \sigma_l'(\xi',\lambda) := \left(\int_{-\infty}^\infty \frac
{\xi_n^{2l}}{\sigma^2(\xi,\lambda)}\,d\xi_n\right)^{-1/2} < \infty\,.\]
Then $D_n^lu(x',0)$ is well-defined as an element of
 $H^{\sigma'_l}(\R^{n-1})$ for
every $u\in H^\sigma(\R^n)$, and there exists a constant $C$, independent
of $u$ and $\lambda$, such that
\begin{equation}\mylabel{3-4}
\| D_n^l u(x',0)\|_{\sigma'_l,\R^{n-1}} \le C \|u\|_{\sigma,\R^n}\,.
\end{equation}
\end{proposition}

We will apply Proposition \ref{3.1} to the case where
$\sigma(\xi,\lambda)$ is given by 
$\Xi_P(\xi,\lambda)$ (see \eqnref{3-2}).

Let one of the functions $\sigma(\xi,\lambda)$ or $\sigma_1(\xi,\lambda)$  for
each $\lambda$ satisfy the condition of Definition \ref{3.0} and
$\sigma(\xi,\lambda) \approx   \sigma_1(\xi,\lambda)$.   The   symbol
$\approx$ means that there exist positive constants $C_1$ and $ C_2$,
independent of $\xi$ and $\lambda$, such that
\[ C_1 \sigma(\xi,\lambda) \le \sigma_1(\xi,\lambda) \le C_2 \sigma(\xi,\lambda)\,.\]

Then the other function also satisfies the condition of Definition
\ref{3.0} and,
evidently, the statement of Proposition \ref{3.1} remains valid,  if
we replace $\sigma$ by the equivalent  function  $\sigma_1$.  In the
following
  we will
construct an equivalent function for $\Xi_P(\xi,\lambda)$ (cf.
\cite{denk-mennicken-volevich}, Section 2). For this purpose we 
introduce some
simple geometric notions connected with the Newton polygon
(see, e.g., \cite{gindikin-volevich}, Chapter 1).

\begin{figure}[ht]
\setlength{\unitlength}{1mm}
\begin{center}
\begin{picture}(60,65)(0,0)
\put(5,10){\vector(1,0){50}}
\put(5,10){\vector(0,1){50}}
\put(52,5){$i$}
\put(0,55){$k$}
\put(5,45){\line(6,-1){15}}
\put(19.8,42.5){\line(2,-1){5}}
\put(24.8,40){\line(1,-1){10}}
\put(34.8,30){\line(1,-3){6.6}}
\put(4.5,9.6){{$\scriptscriptstyle \bullet$}}
\put(4.4,44.6){{$\scriptscriptstyle \bullet$}}
\put(19.2,42){{$\scriptscriptstyle \bullet$}}
\put(24.2,39.5){{$\scriptscriptstyle \bullet$}}
\put(34.2,29.5){{$\scriptscriptstyle \bullet$}}
\put(40.9,9.6){{$\scriptscriptstyle \bullet$}}
\put(9,46){$\Gamma_1$}
\put(21.4,43){$\Gamma_2$}
\put(29,38){$\ddots$}
\put(40,22){$\Gamma_S$}
\put(5,2){\parbox{50mm}{\begin{center}
{\small Fig. 1}\end{center}}}
\end{picture}
\end{center}
\end{figure}

Let $\Gamma_1,\ldots,\Gamma_S$ be the sides of the Newton polygon 
not lying on the coordinate axes and indexed in the clockwise direction
(cf. Fig. 1). Suppose that
\[ (0,0),\, (a_1,b_1), \ldots, (a_{S+1},b_{S+1})\,,
\quad a_1=0,\quad b_{S+1}=0\,, \]
are the vertices of the polygon $N(P)$. Then  the side
$\Gamma_s$ is given by
\begin{equation}\mylabel{3-5a}
\Gamma_s =
\{ (a,b)\in\R^2: 1\cdot a+ r_s\cdot b= d_s\}\quad (s=1,\ldots,S)
\end{equation}
where $ r_s=(a_{s+1}-a_s)/(b_s-b_{s+1})$.
The vector $(1,r_s)$ is an exterior normal to the side
$\Gamma_s$, where we admit $r_1=\infty$ if $\Gamma_1$ is horizontal.
Further we have $r_S=0$  in  the case that $\Gamma_S$ is vertical.  In what
follows we will  suppose  that  $\Gamma_S$  is  not  vertical.
Since $N(P)$ is convex, we have
\[     \infty \ge r_1 > \ldots >r_S > 0\,.\]
The $r_s$-principal part of $P$ is defined by
\begin{equation}\mylabel{3-5'}
P_{r_s}(\xi,\lambda) := \sum_{|\alpha|+r_sk = d_s} a_{\alpha k}
\xi^\alpha\lambda^k\,.
\end{equation}
Here $d_s$ is the so-called $r_s$-degree of $P$ which may be defined
by
\begin{equation}\mylabel{3-5b}
d_s := \max_{(a,b)\in N(P)} (1\cdot a+r_s\cdot b)\,.
\end{equation}

Now we set
\[     \Xi_{(s)}(\xi,\lambda)= |\xi|^{-a_s}\;|\lambda|^{-b_{s+1}}
     \sum_{i+r_sk=d_s} |\xi|^i|\lambda|^k\,. \]
This function will be a polynomial of $|\xi|$ and $|\lambda|$.

Repeating the argument in \cite{gindikin-volevich},  Theorem  1.1.3,
we can prove that
\begin{equation}\mylabel{3-5c}
    \prod_{s=1}^S \Xi_{(s)}(\xi,\lambda)= \sum_{s=1}^S
      |\xi|^{a_s}|\lambda|^{b_s} + \dots\;,
\end{equation}
where the dots denote the sum of monomials $|\xi|^i|\lambda|^k$
with $(i,k) \in N(P)$. For $|\lambda| \ge1$ the right-hand side
can be estimated from below by
\[    1+ \sum_{s=1}^S
|\xi|^{a_s}|\lambda|^{b_s} \, .\]
This function    can    be    estimated    from    below    by
$\Xi_P(\xi,\lambda)$ (see  \cite{denk-mennicken-volevich},
Subsection 3.2). From
this it follows that the left-hand side of \eqnref{3-5c} is equivalent to
$\Xi_P$. Denote by $2m_s$ the largest degree of $|\xi|$ in
$\Xi_{(s)}$. It is obvious that $\Xi_{(s)}$ is equivalent to
$ (|\xi|+|\lambda|^ \frac{1}{r_s})^{2m_s}$, and consequently
\begin{equation}\mylabel{3-5}
\Xi_P(\xi,\lambda) \approx \prod_{s=1}^S \left(
|\xi|^2 + |\lambda|^{\frac 2 {r_s}}\right)^{m_s}\,.
\end{equation}

We will   suppose   further,  as  in  the  case  of  parabolic
polynomials (cf.\ \cite{gindikin-volevich},  Chapter  2), that   $m_1,
\dots, m_S$ are integers. 

\begin{remark}\mylabel{3.0aa}{\rm
 In the case $r_1=\infty$ (i.e. $\Gamma_1$ is horizontal)
\eqnref{3-5'} and \eqnref{3-5b}
 have no sense and \eqnref{3-5'} should be replaced by
\[ P_{r_1}:=\sum_{|\alpha|=a_2} a_{\alpha b_1}\xi^{\alpha}\lambda^
{b_1}\,.\]
As for the equivalence \eqnref{3-5}, it will be valid for $|\lambda|
>\lambda_0$ with arbitrary $\lambda_0>0$ and the equivalence constants,
of course, depend on $\lambda_0$.}
\end{remark}

\begin{remark}\mylabel{3.0a} {\rm The fact that $\Xi(\xi,\lambda)$
satisfies the condition of Definition \ref{3.0} is an immediate
consequence of \eqnref{3-5} as this condition is fulfilled for each
factor on the right-hand side.}
\end{remark}

\begin{remark}\mylabel{3.1a}{\rm From  \eqnref{3-5} it follows that
the $r_s$-degree $d_s$ (cf. \eqnref{3-5b}) is given by
\begin{equation}\mylabel{3-5e}
d_s = 2\left(\sum_{j=1}^s m_j + \sum_{j=s+1}^S \frac{r_s}{r_j}m_s
\right) \,.
\end{equation}
To see this, we use the relation
\begin{equation}\mylabel{3-5f}
\Xi_P(t\xi, t^{r_s}\lambda) = t^{d_s} \Xi_{P_{r_s}}(\xi,\lambda)
+ o(t^{d_s}),\quad t\to+\infty\,,
\end{equation}
cf. \cite{gindikin-volevich}, Section 1.1.2. In our case we obtain, denoting
the right-hand side of \eqnref{3-5e} by $d_s'$,
\begin{eqnarray*}
\lefteqn{\Xi_P(t\xi,t^{r_s}\lambda)  = \prod_{j=1}^S \left(
t^2 |\xi|^2 + t^{2\frac{r_s}{r_j}} |\lambda|^{\frac 2{r_j}}\right)^{m_j}}\\
& = & t^{d_s'} \prod_{j=1}^s  \left(
|\xi|^2 + t^{2(\frac{r_s}{r_j}-1)} |\lambda|^{\frac 2{r_j}}\right)^{m_j}
 \prod_{j=s+1}^S  \!\!\left( t^{2(1-\frac{r_s}{r_j})}
|\xi|^2 +  |\lambda|^{\frac 2{r_j}}\right)^{m_j}\\
& = & t^{d_s'} \Xi_{P_{r_s}}(\xi,\lambda) + o(t^{d_s'})\,,
\end{eqnarray*}
which shows $d_s=d_s'$.}
\end{remark}

Now we will describe the trace spaces of the spaces $H^\Xi$. For this we
use the following lemma:

\begin{lemma}\mylabel{3.2}
Let $1\le a_1<a_2<\ldots<a_S<\infty$ and $m_1,\ldots,m_S\in\N$. For 
$l\in\N$ with $0\le l < 2(m_1+\ldots+m_S)$ define the index $\kappa$
by
\begin{equation}\mylabel{3-6}
2m_1+\ldots + 2m_{\kappa-1} \le l < 2m_1+\ldots+2m_\kappa\,.
\end{equation}
Then there exists a constant $C>0$, independent of $a_1,\ldots,a_S$, such
that
\begin{equation}\mylabel{3-7}\begin{array}{l}{\displaystyle{ 
C^{-1} a_\kappa^{2l+1-4m_1-\ldots-4m_\kappa}\prod_{s=\kappa+1}^S
a_s^{-4m_s} \le \int_{-\infty}^\infty \frac{t^{2l}}{\prod_{s=1}^S
(t^2+a_s^2)^{2m_s}}\,dt}}\\ {\displaystyle{\hspace*{4cm}
 \le  C a_\kappa^{2l+1-4m_1-\ldots-4m_\kappa}
\prod_{s=\kappa+1}^Sa_s^{-4m_s}\,.}}\end{array}
\end{equation}
In the case $0\le l<2m_1$, we set $m_0=0$ in \eqnref{3-6}.
\end{lemma}

\begin{proof}
Substituting in the integral $t=a_S\tau$, we obtain
\begin{eqnarray*}
 I &:=& \int_{-\infty}^\infty t^{2l} \prod_{s=1}^S(t^2+a_s^2)^{-2m_s}
dt \\
&=& 2 a_S^{2l+1-4m_1-\ldots-4m_S}\int_0^\infty t^{2l} \prod_{s=1}^S
\left( t^2 + \Big(\frac{a_s}{a_S}\Big)^2\right)^{-2m_s}dt.
\end{eqnarray*}
For $t\ge 1$ we use
\[ t^{2l}(1+t^2)^{-2 m_1-\ldots -2m_S} \le t^{2l} \prod_{s=1}^S 
\left( t^2 + \left(\frac{a_s}{a_S}\right)^2\right)^{-2m_s}\le 
t^{2l-4m_1-\ldots-4m_S}\,.\]
As $l<2\sum_{s=1}^S m_s$, the left-hand and right-hand side of this 
inequality are integrable functions over $[1,\infty)$, and we obtain
\[ C_1^{-1}\le \int_1^\infty t^{2l} \prod_{s=1}^S 
\left( t^2 + \left(\frac{a_s}{a_S}\right)^2\right)^{-2m_s}\!dt\le C_1\]
for some $C_1>0$. 

For $0\le t\le 1$ we have $1\le 1+t^2 \le 2$, and therefore
\[ \int_0^1\ldots dt \approx \int_0^1 t^{2l} \prod_{s=1}^{S-1} \left(
t^2 + \frac{a_s^2}{a_S^2}\right)^{-2m_s}\, dt\,.\]
Now we substitute $t=\frac{a_{S-1}}{a_S}\tau$ and get
\[ \int_0^1 \!\!\ldots dt \approx
\left(\frac{a_{S-1}}{a_S}\right)^{\!\!
2l+1-4m_1
-\ldots-4m_{S-1}}\!\!\int_0^{\frac{a_S}{a_{S-1}}} \!\! t^{2l} 
\prod_{s=1}^{S-1} \left(t^2 +
  \frac{a_s^2}{a_{S-1}^2}\right)^{\!\!-2m_s}\!\!
 dt\,.\]
Again we split up $\int_0^{\frac{a_S}{a_{S-1}}}\ldots = \int_0^1\ldots
+ \int_1^{\frac{a_S}{a_{S-1}}}\ldots$ and use an estimate of the form
$C_2^{-1} \le \int_1^{\frac{a_S}{a_{S-1}}}\ldots\le C_2$ for the second
integral.

Proceding in this way, we receive
\begin{eqnarray*}
 I& \approx&  a_S^{2l+1-4m_1-\ldots-4m_S}\left(\frac{a_{S-1}}{a_S}\right)
^{2l+1-4m_1-\ldots-4m_{S-1}}\cdot\ldots\\
& & \cdot \left(\frac{a_{\kappa
}}{a_{\kappa+1}}\right)^{2l+1-4m_1-\ldots-4m_\kappa}\int_0^{\frac{a_{\kappa
+1}}{a_\kappa}}   t^{2l} 
\prod_{s=1}^\kappa \left(t^2 + \frac{a_s^2}{a_\kappa^2}
\right)^{-2m_s}\, dt\,.
\end{eqnarray*}
For the last integral we use
\begin{eqnarray*}
t^{2l} (t^2+1)^{-2m_1-\ldots-2m_\kappa} &\le& t^{2l} \prod_{s=1}^\kappa
 \left(t^2 + \frac{a_s^2}{a_\kappa^2}\right)^{-2m_s}\\
&\le& t^{2l-4m_1-\ldots-4m_{\kappa-1}}(t^2+1)^{-2m_\kappa}\,.
\end{eqnarray*}
As $2m_1+\ldots+2m_{\kappa-1}\le l < 2m_1+\ldots+2m_\kappa$, the
left-hand and the right-hand side of this inequality are integrable 
functions on $[0,\infty)$. Therefore
\[ I \approx a_\kappa^{2l+1-4m_1-\ldots-4m_\kappa} a_{\kappa+1}^{-4
m_{\kappa+1}}\cdot \ldots \cdot a_S^{-4m_S}\,.\]
\end{proof}

\begin{remark}\mylabel{3.2a}{\rm Using the substitution $t=a_1\tau$, it 
    is easily seen that the condition $a_1\ge 1$ in Lemma \ref{3.2}
    may be replaced by $a_1>0$.}
\end{remark}

As in the Introduction, we denote by $\Xi^{(-l)}_P(\xi,\lambda)$ the
function corresponding to the Newton polygon which is constructed from 
$N(P)$ by a shift of length $l$ to the left parallel to the
abscissa. We preserve the notation $H^{\Xi_P^{(-l)}}(\R^{n-1})$ for the
spaces in $\R^{n-1}$ corresponding to the weight functions
$\Xi_P^{(-l)}(\xi', \lambda): = \Xi_P^{(-l)}(\xi',0,\lambda)$.

\begin{lemma}\mylabel{3.2b} Let $\lambda_0>0$. Then for $|\lambda| \ge 
  \lambda_0$ we have
\begin{equation}\mylabel{3-7a}
\sigma_l'(\xi',\lambda)\approx \Xi^{(-l-\frac 1 2)}(\xi',\lambda)\,,
\end{equation}
where $\sigma_l'$ is defined by
\begin{equation}\mylabel{3-7b}
\sigma_l'(\xi',\lambda) := \left( \int_{-\infty}^\infty \frac{\xi_n^{2l}}{
\Xi_P^2(\xi,\lambda)}\,d\xi_n\right)^{-\frac 1 2}\,.
\end{equation}
\end{lemma}

\begin{proof} Instead of $\Xi_P$ we use the right-hand side of
  \eqnref{3-5}. 
From Lemma \ref{3.2} with $a_s^2 = |\xi'|^2+
|\lambda|^{\frac{2}{r_s}}$ we obtain (see Remark \ref{3.2a}) that
\begin{equation}\mylabel{3-10}
\sigma_l'(\xi',\lambda) \approx \left( |\xi'|^2 + |\lambda|^{\frac 2{r_\kappa}}
\right)^{m_1+\ldots+m_\kappa-\frac l 2-\frac 1 4}\prod_{s=\kappa +1}^S
\left( |\xi'|^2 + |\lambda|^{\frac 2{r_s}}\right)^{m_s}\,,
\end{equation}
where $\kappa$ is chosen according to Lemma \ref{3.2}. From Remark
\ref{3.1a} applied to $\sigma_l'(\xi',\lambda)$ we see that the
sides of the Newton polygon corresponding to the weight function
\eqnref{3-10} are given by
\[ \Gamma_j = \{ (a,b)\in\R^2: a+r_j b = d_j'\} \]
with $d_j' = d_j-l-\frac l 2\; (j=\kappa,\ldots,S)$. But this means
that the Newton polygon for $\sigma_l'$ is constructed from $N(P)$ by
a shift of $l+\frac 1 2$ to the left, i.e. we have
$\sigma_l'(\xi',\lambda)\approx \Xi_P^{(-l-\frac 1 2)}(\xi',
\lambda)$.
\end{proof}

The following theorem is an immediate consequence of Proposition
\ref{3.1} and Lemma \ref{3.2b}.

\begin{theorem}\mylabel{3.3} For every  $\lambda_0>0$ there exists a
  constant $C>0$, independent of $u$ and $\lambda$, such that
\begin{equation}\mylabel{3-8}
\| D_n^l u(x',0)\|_{\Xi_P^{(-l-\frac 1 2)},\R^{n-1}}
 \le C \| u \|_{\Xi_P,\R^n}\quad(l=0,\ldots,2m_1+\ldots+2m_S-1)
\end{equation}
holds for $u\in H^{\Xi_P}(\R^n)$ and  $\lambda\in\C$ with
$|\lambda|\ge\lambda_0$.
\end{theorem}

In the following, we will also consider the function spaces in the
half space $\R^n_+$ which
correspond to Newton polygons. Using the binomial formula, it is
easily seen that 
\begin{equation}\mylabel{3-9}
\Xi^2_P(\xi,\lambda) \approx \sum_{l=0}^{M} \xi_n^{2l}\,\big(\Xi_P^{(-l)}
(\xi',\lambda)\big)^2
\end{equation}
where $M=2m_1+\dots +2m_S$. From this it  follows that we can take
\begin{equation}\mylabel{3-10a}
\Big(\sum_{l=0}^M \int_{-\infty}^{\infty}
\|(D_n^l u)(\cdot,x_n)\|^2_{\Xi^{(-l)}_
{P}, \R^{n-1}}\, dx_n\Big)^{1/2}
\end{equation}
as an equivalent norm in $H^{\Xi_P}(\R^n)$. Replacing the integral over $\R$
by the integral over ${x_n \ge 0}$ we define a norm in
$H^{\Xi_P}(\R^n_+)$.

To define the space $H^{\frac 1 {\Xi_P}}(\R^n_+)$, we use the more general
approach which can be found, e.g., in \cite{volevich-panejah}. Let
$\sigma(\xi)$ be a weight function fulfilling the condition in Definition 
\ref{3.0}. Denote by $H^\sigma(\R^n)_{\pm}$ the subspace of $H^\sigma(\R^n)$
consisting of elements with supports in the closure of
$\R^n_{\pm}$. Then we define
\begin{equation}\mylabel{3-11}
 H^\sigma(\R^n_+) = H^\sigma(\R^n) / H^\sigma(\R^n)_- 
\end{equation}
endowed with the natural quotient norm
\begin{equation}\mylabel{3-12}
 \| f\|_{\sigma,\R^n_+} = \inf_{f_-\in H^\sigma(\R^n)_-} \| f_0 +
f_-\|_{\sigma,\R^n}\,,
\end{equation}
where $f_0$ is an arbitrary representative of the conjugacy class of
$f$.

Suppose that $\sigma(\xi',\xi_n)$ for a fixed $\xi'\in\R^{n-1}$ can be
extended as a holomorphic function in $\xi_n$ of polynomial growth in
the lower half-plane $\Im \xi_n<0$. In this case the quotient norm of
$f\in H^\sigma(\R^n_+)$ coincides with the norm
\begin{equation}\mylabel{3-13}
 \|\sigma(D', D_n)f_0\|_{L_2(\R^n_+)} 
\end{equation}
which does not depend on the choice of the element $f_0$ in the
conjugacy class. In \eqnref{3-13} the pseudo-differential operator (ps.d.o.)
$\sigma(D', D_n)=\sigma(D)$ is defined by 
\[ \sigma(D) f := F^{-1} \sigma(\xi) (Ff)(\xi)\,\]
In the case when 
\[ \sigma\approx \prod_{j=1}^S (|\xi|^2 + |\lambda|^{2/r_j})^{m_j} \]
we replace $\sigma$ in the definition of $H^\sigma(\R^n_+)$ by
\[ \prod_{j=1}^S \Big(i\xi_n + (|\xi|^2 + |\lambda|^{2/r_j})^{1/2}
\Big)^{2m_j}\,.\]

\mysection{3}{The zeros of the symbol}

Now we come back to the operator pencil \eqnref{1-1} and consider
the corresponding model problem with constant coefficients and without 
lower order terms.  Let $A(\xi,\lambda)$
be a polynomial in $\xi\in\R^n$ and $\lambda\in\C$ of the form
\begin{equation}\mylabel{2-1a}
A(\xi,\lambda) = A_{2m}(\xi) + \lambda A_{2m-1}(\xi) +\ldots
+\lambda^{2m-2\mu} A_{2\mu}(\xi)\,,
\end{equation}
where $A_j(\xi)$ is a homogeneous polynomial in $\xi$ of degree $j$.

\begin{definition}\mylabel{2.0}{\rm  The polynomial $A(\xi,\lambda)$ is
called $N$-elliptic with parameter in $[0,\infty)$ if the estimate
\begin{equation}\mylabel{2-1}
|A(\xi,\lambda)|\ge C\,|\xi|^{2\mu}
\;(\lambda + |\xi|)^{2m-2\mu}\quad ( \xi\in
\R^n, \lambda\in [0,\infty))
\end{equation}
holds with a constant $C$ independent of $\xi$ and $\lambda$.}
\end{definition}

\begin{lemma}\mylabel{2.1} The polynomial $A(\xi,\lambda)$ is
$N$-elliptic with parameter in\newline  $[0,\infty)$ if   and only
if the following conditions are satisfied:\\[2mm]
{\rm (i)} $A_{2m}(\xi)$ is elliptic, i.e. $A_{2m}(\xi)\not= 0 $ for $\xi
\in \R^n\backslash\{0\}$.\\[2mm]
{\rm (ii)} $A_{2\mu}(\xi)$ is elliptic.\\[2mm]
{\rm (iii)} $A(\xi,\lambda)\not=0$ for $\xi\in\R^n\backslash\{0\}$ and 
$\lambda\in[0,\infty)$.
\end{lemma}

\begin{proof} From (\ref{2-1}) we trivially obtain condition (iii) and,
setting $\lambda=0$, condition (i). Taking $\epsilon=\frac 1 \lambda$ and
dividing (\ref{2-1}) by $\epsilon^{2\mu-2m}$, we receive
\begin{equation}\mylabel{2-2}
|A_{2\mu}(\xi) + \epsilon A_{2\mu+1}(\xi) + \ldots + \epsilon^{2m-2\mu}
A_{2m}(\xi)|\ge C |\xi|^{2\mu} (1+\epsilon|\xi|)^{2m-2\mu}\,.
\end{equation}
Taking the limit for $\epsilon\to 0$, we obtain (ii).

Now let conditions (i)--(iii) be fulfilled. For $\xi\in\R^n\backslash
\{0\}$ we write $A(\xi,\lambda)$
in the form
\begin{equation}\mylabel{2-3}
A(\xi,\lambda) = A_{2\mu}(\xi) B_{2m-2\mu}(\xi,\lambda)
\end{equation}
with
\begin{equation}\mylabel{2-4}
B_{2m-2\mu}(\xi,\lambda) = \frac{A_{2m}(\xi)}{A_{2\mu}(\xi)} +
\lambda \frac{A_{2m-1}(\xi)}{A_{2\mu}(\xi)} + \ldots + \lambda^{2m-2\mu}\,.
\end{equation}
The coefficients of $B_{2m-2\mu}(\xi,\lambda)$ (considered as a
polynomial in $\lambda$) are homogeneous functions in $\xi\in
\R^n\backslash\{0\}$, and therefore $B(\xi,\lambda)$ is a homogeneous 
function in $(\xi,\lambda)$ of degree $2m-2\mu$. From this and from
conditions (ii) and (iii) it follows that
\begin{equation}\mylabel{2-5}
|A_{2\mu}(\xi) | \ge C |\xi|^{2\mu}\,,\quad |B_{2m-2\mu}(\xi,\lambda)|
\ge C(\lambda + |\xi|)^{2m-2\mu}\,.
\end{equation}
Multiplying these estimates, we see that $A$ is $N$-elliptic with parameter
in $[0,\infty)$.
\end{proof}

Denote by $\tau_j(\xi',\lambda)$ ($j=1,\ldots,2m$)
 the zeros of the algebraic
equation
\begin{equation}\mylabel{2-6}
A(\xi',\tau,\lambda)=0 \quad \big(\xi'\in\R^{n-1}\backslash \{0\}, 
\lambda\in [0,\infty)\big)\,.
\end{equation}

Due to Lemma \ref{2.1} (iii), this equation has no real roots. The number
$m_+$ of roots with positive imaginary part is independent of $(\xi',
\lambda)$ and therefore coincides with the corresponding number for
$\lambda=0$. It is easily seen (cf. \cite{agranovich96}, Section 1.2) that in
the case $n>2$ the set $\{(\xi',\lambda): \xi'\in\R^{n-1}\backslash
\{0 \}, \;
\lambda\in [0,\infty)\}$ is connected, and therefore
 we have $m_+=m$. In the case $n\le 2$ the relation $m_+=m$
is an additional condition which will be assumed to hold in the 
following. We denote the roots of $A(\xi',\tau,\lambda)$ with positive
imaginary part by $\tau_1(\xi',\lambda),\ldots,\tau_m(\xi',\lambda)$.

To investigate the elliptic pencil corresponding to 
 $A(\xi',\tau,\lambda)$ we will need
an additional assumption which is closely related to the condition of
regularity of degeneration in the theory of singular perturbations (cf.
Vishik-Lyusternik \cite{vishik-lyusternik}, Section 1.1). 
To formulate this assumption we consider
the auxiliary polynomial of degree $2m-2\mu$ given by
\begin{equation}\mylabel{2-7}
Q(\tau) := \tau^{-2\mu} A(0,\tau,1)\,.
\end{equation}
From  inequality (\ref{2-1}) with $\xi'=0$ and $\lambda=1$ we 
obtain for $\tau\not=0$ the estimate
\begin{equation}\mylabel{2-8}
|Q(\tau)| \ge C(|\tau|+1)^{2m-2\mu}
\end{equation}
with a constant independent of $\tau$. By continuity we obtain that
$Q(0)\not=0$, and thus $Q(\tau)$ has no real roots. 

\begin{definition}\mylabel{2.1a}{\rm The polynomial $A(\xi',\tau,\lambda)$
is said to degenerate regularly for $\lambda\to\infty$ if the 
polynomial $Q(\tau)$ defined in \eqnref{2-8} has exactly $m-\mu$ roots
with positive imaginary part (counted according to their multiplicities).}
\end{definition}

\begin{remark}\mylabel{2.2} {\rm a) Suppose that the polynomial $A(\xi
,\lambda)$ contains only terms of even order, i.e.
\begin{eqnarray}
A(\xi,\lambda)& =& A_{2m}(\xi)+\lambda^2 A_{2m-2}(\xi) +\ldots\nonumber\\
& & + \; \lambda^{2m-2\mu-2} A_{2\mu+2}(\xi) + \lambda^{2m-2\mu}A_{2\mu}
(\xi)\,.\mylabel{2-10}
\end{eqnarray}
Then the polynomial $Q(\tau)$ is a polynomial of degree $m-\mu$ in the
variable $\tau^2$ and $A(\xi,\lambda)$ degenerates regularly for 
$\lambda\to\infty$.

b) (Cf. \cite{vishik-lyusternik}, Lemma 3.4.) Assume that
$A(\xi,\lambda)$ is the symbol of a differential operator $\tilde A
(\frac\partial{\partial x_1},\ldots,\frac\partial{\partial x_n},
\lambda)$ with real coefficients. Then the polynomials of even order
$A_{2m-2j}(\xi)\;(j=0,\ldots,m-\mu)$ are real and the polynomials
of odd order $A_{2m-2j-1}(\xi)\;(j=0,\ldots,m-\mu-1)$ are purely
imaginary. Assume that $\tilde A$ is strongly elliptic, i.e. we have
\begin{equation}\mylabel{2-11}
\Re A(\xi,\lambda) \ge C |\xi|^{2\mu} (\lambda + |\xi|)^{2m-2\mu}\,.
\end{equation}
Then we obtain that $\Re A = A_{2m} + \lambda^2 A_{2m-2} + \ldots
+ \lambda^{2m-2\mu}A_{2\mu}$ satisfies (\ref{2-1}), 
and due to part a) the polynomial 
$\Re Q(\tau)$ has $m-\mu$ roots with positive imaginary part and $m-\mu$
roots with negative imaginary part. Since the polynomial
\begin{equation}\mylabel{2-12}
Q_\delta(\tau) := \Re Q(\tau) + \delta i \Im Q(\tau)\quad (0\le\delta
\le 1)
\end{equation}
satisfies
\begin{equation}\mylabel{2-13}
\Re Q_\delta(\tau) \ge C (|\tau|+1)^{2m-2\mu}\quad  (0\le\delta
\le 1)\,,
\end{equation}
the number of roots of $Q_\delta$ in the upper half complex plane does
not depend on $\delta\in[0,1]$, and $A(\xi,\lambda)$ degenerates 
regularly for $\lambda\to\infty$.}
\end{remark}

\begin{lemma}\mylabel{2.3}  Let the polynomial $A(\xi,\lambda)$ in
\eqnref{2-1a} be $N$-elliptic with parameter in $[0,\infty)$ and
assume that $A$ degenerates regularly for $\lambda\to\infty$.
 Then, with a suitable numbering of the roots $\tau_j(\xi',
\lambda)$ of $A(\xi',\tau,\lambda)$ with positive imaginary part, we 
have:\\
{\rm (i)} Let $S(\xi')= \{ \tau_1^0(\xi'),\ldots,\tau_\mu^0(\xi')\}$ be the 
set of all zeros of $A_{2\mu}(\xi',\tau)$ with positive imaginary
part. Then for all $r>0$ there exists a $\lambda_0>0$ such that the
distance between the sets $\{
\tau_1(\xi',\lambda),\ldots,\tau_\mu(\xi',\lambda)\}$ and $S(\xi')$ is 
less than $r$ for all $\xi'$ with $|\xi'|=1$ and all
$\lambda\ge\lambda_0$.\\
{\rm (ii)}  Let $\tau_{\mu+1}^1,\ldots,\tau_m^1$ be the roots of the polynomial
$Q(\tau)$ (cf. \eqnref{2-7})
 with positive imaginary part. Then
\begin{equation}\mylabel{2-16}
\tau_j(\xi',\lambda) = \lambda\tau_j^1 + \tilde\tau_j^1(\xi',\lambda)
\quad(j=\mu+1,\ldots,m)\,,
\end{equation}
and there exist  constants $K_j$ and $\lambda_1$, independent of
$\xi'$ and $\lambda$,
such that for $\lambda\ge \lambda_1$ the inequality
\begin{equation}\mylabel{2-17}
|\tilde\tau_j^1(\xi',\lambda)|\le K_j |\xi'|^{\frac 1{k_1}}\,
\lambda^{1-\frac 1{k_1}}\quad  (|\xi'|\le\lambda)
\end{equation}
holds, where $k_1$ is the maximal multiplicity of the roots of $Q(\tau)$.
\end{lemma}

\begin{proof}
(i) We write $\xi'=\rho\omega$ with $|\omega|=1$ and set
$\tilde \tau=\frac\tau\rho,\; \epsilon=\frac\rho\lambda$. After division
of $A(\xi',\tau,\lambda)$ by $\lambda^{2m-2\mu}\rho^{2\mu}$ we obtain 
the equation
\begin{equation}\mylabel{2-18}
B(\omega,\tilde\tau,\epsilon) := A_{2\mu}(\omega,\tilde\tau)
+\epsilon A_{2\mu+1}(\omega,\tilde\tau) + \ldots + \epsilon^{2m-2\mu}
A_{2m}(\omega,\tilde\tau)=0\,.
\end{equation}
First we fix $\omega$ with $|\omega|=1$. Let
$\tilde\tau_j=\ldots=\tilde \tau_{j+p-1}$ be a zero of
$B(\omega,\tilde\tau,0)=A_{2\mu}(\omega, \tilde \tau)$ of multiplicity 
$p$. Then there exists an $\alpha=\alpha(\omega)>0$ such that 
\begin{equation}\mylabel{2-19}
\frac{1}{2\pi i} \int_{|z-\tilde\tau_j|=\alpha} \frac{\frac{d}{dz} 
B(\omega,z,\epsilon)}{B(\omega,z,\epsilon)} \,dz =
\frac{1}{2\pi i} \int_{|z-\tilde\tau_j|=\alpha} \frac{\frac{d}{dz} 
B(\omega,z,0)}{B(\omega,z,0)} \,dz = p
\end{equation}
holds for all $\epsilon<\epsilon_0=\epsilon_0(\omega)$. Therefore, for 
every $\epsilon<\epsilon_0$ the equation \eqnref{2-18} has exactly $p$ 
roots in $\{z\in\C: |z-\tilde\tau_j|<\alpha\}$ which we denote by
$\tilde\tau_j(\omega,\epsilon), \ldots,
\tilde\tau_{j+p-1}(\omega,\epsilon)$. Proceding in this way for all
zeros of $A_{2\mu}(\omega,\tilde\tau)$, we obtain the set
$S(\omega,\epsilon) := \{\tilde\tau_1(\omega,\epsilon), \ldots,
\tilde\tau_\mu(\omega,\epsilon)\}$ of zeros of
$B(\omega,\tilde\tau,\epsilon)$.

Now we assume that the statement in (i) is false. Then there exists a
sequence $(\omega_n)_{n\ge 1}$ with $|\omega_n|=1$ and a constant
$C>0$ such that ${\rm dist}(S(\omega_n), S(\omega_n,\epsilon_n)) \ge
C$ for all $n\ge 1$ where we have set $\epsilon_n=\frac 1 n$.
Due to compactness, we may assume that $\omega_n$ converges to
$\omega_0$. As the zeros of $A_{2\mu}(\omega,\tilde\tau)$ depend
continuously on $\omega$, we obtain for large $n$ that
\begin{equation}\mylabel{2-20}
{\rm dist}(S(\omega_0), S(\omega_n,\epsilon_n)) \ge\frac C2\,.
\end{equation}
But from the same considerations as above we see that for every
sufficiently small $\alpha>0$ there exists an $\epsilon_0 =
\epsilon_0(\omega_0)$ and an $s>0$ such that
$B(\omega,\tilde\tau,\epsilon)$ has exactly $\mu$ roots in $\bigcup_j
\{z\in\C: |z-\tilde\tau_j(\omega_0)|<\alpha\}$ for all $|\omega-\omega_0| < 
s$ and $0<\epsilon<\epsilon_0$. Taking $\alpha<\frac C 2$, we obtain a 
contradiction to \eqnref{2-20}.

(ii)  We set $\tilde\tau=\frac\tau\lambda$ and $\epsilon=\frac{|\xi'|}\lambda$
and obtain the equation $B(\omega,\tilde\tau,\epsilon) :=
A(\epsilon\omega, \tilde\tau,1) = 0$ with $\omega:=
\frac{\xi'}{\epsilon}$. First we fix $\omega$ with $|\omega|=1$. We write
\begin{equation}\mylabel{2-21}
0 = B(\omega,\tilde\tau,\epsilon) = A(0,\tilde\tau,1) + \sum_{k=1}^{2m}
 \Big(\frac\partial{\partial\epsilon}\Big)^k B(\omega,\tilde\tau,0)
\;\frac{\epsilon^k}{k!}\,.
\end{equation}
Let $\tau^1_j=\ldots=\tau^1_{j+p-1}$ be a zero of $Q(\tau)$ of multiplicity
$p$. Then we know from the theory of algebraic functions that there
exist $p$ roots $\tilde\tau_j(\omega,\epsilon),\ldots,\tilde
\tau_{j+p-1}(\omega,\epsilon)$ of $B(\omega,\tilde\tau,\epsilon)$ for which 
 we have an expansion (Puiseux series) of the form
\begin{equation}\mylabel{2-22}
\tilde\tau_s(\omega,\epsilon) = \tau_j^1 + \sum_{k=1}^\infty
c_{jk}(\omega) \epsilon^{k/p}\quad (s=j,\ldots,j+p-1)
\end{equation}
(cf., e.g., \cite{fuks-levin}, Section 7). In formula \eqnref{2-22} we have 
to take the $p$ different branches of the function $\epsilon^{\frac 1 p}$
to obtain the zeros $\tilde\tau_j(\epsilon), \ldots, 
\tilde\tau_{j+p-1}(\epsilon)$. The series on the right-hand
side is a holomorphic function in $\epsilon^{\frac 1 p}$ for 
$|\epsilon|\le\epsilon_1(\omega)$ for some  $\epsilon_1(\omega)>0$. 

From the construction of the Puiseux series (cf. \cite{fuks-levin},
Section 8) we know that the coefficients $c_{jk}(\omega)$ in the
series \eqnref{2-22} depend continuously on the coefficients of the
polynomial $B(\omega,\tilde\tau,\epsilon)$ and therefore on
$\omega$. Thus there exists an $\epsilon_1 >0$, independent of
$\omega$, such that the right-hand side of \eqnref{2-22} is a
holomorphic function in $\epsilon^{\frac 1 p}$ for $|\epsilon|\le
\epsilon_1$. As the function 
\[ (\tilde\tau_j(\omega,\epsilon)
-\tau_j^1) \epsilon^{-\frac 1 p} = \sum_{k=1}^\infty c_{jk}(\omega)
\epsilon^{\frac{k-1}{p}} \]
is continuous in $\omega$ and $\epsilon$ for $|\omega|=1$ and $0\le
\epsilon \le \epsilon_0$, it is bounded by some constant $K_1$,
independent of $\omega$ and $\epsilon$, which finishes the proof of
part (ii).
\end{proof}

\mysection{4}{Estimates for ordinary differential equations}

In this section we consider the polynomial $A(\xi,\lambda)$
given by \eqnref{2-1a} and assume that this polynomial is $N$-elliptic with
parameter in $[0,\infty)$ and degenerates regularly for $\lambda\to\infty$.
The Newton polygon corresponding
to $A$ has the shape indicated in Figure 2 with $r=2m$ and
$s=2\mu$.

\begin{figure}[ht]
\setlength{\unitlength}{1mm}
\begin{center}
\begin{picture}(90,70)(-5,5)
\put(5,15){\vector(1,0){80}}
\put(5,15){\vector(0,1){50}}
\put(82,10){$i$}
\put(0,60){$k$}
\put(5,41){\line(1,0){35}}
\put(40,41){\line(1,-1){25.8}}
\put(4.4,40.6){{$\scriptscriptstyle \bullet$}}
\put(39.4,40.6){{$\scriptscriptstyle \bullet$}}
\put(65,14.6){{$\scriptscriptstyle \bullet$}}
\put(40,14){\line(0,1){2}}
\put(65,10){$r$}
\put(39,10){$s$}
\put(-6,40){$r-s$}
\put(5,2){\parbox{60mm}{\begin{center}
{\small Fig. 2. The Newton polygon $N_{r,s}$.}\end{center}}}
\end{picture}
\end{center}
\end{figure}

For fixed $\xi'\in\R^{n-1}$, $\lambda\in[0,\infty)$ and $j=1,\ldots
,m$ we consider
 the ordinary differential
equation on the half-line
\begin{eqnarray}
A(\xi', D_t,\lambda)\, w_j(t)  & = &0 \quad\quad (t>0)\,,
\mylabel{4-2}\\
D_t^{k-1}\, w_j(t)|_{t=0} & = & \delta_{jk}  \quad(k=1,\ldots,m)\,,
\mylabel{4-3}\\
w_j(t) & \to & 0 \quad\quad (t \to + \infty)\,.\nonumber
\end{eqnarray}
Here $D_t$ stands for $-i\frac\partial{\partial t}$.

\begin{theorem}\mylabel{4.1} For every $\xi'\in\R^{n-1}\backslash\{0\}$
 and $\lambda\in [0,\infty)$ the ordinary differential equation 
\eqnref{4-2}--\eqnref{4-3} has a unique solution $w_j(\xi',t,\lambda)$,
and the estimate
\begin{equation}\mylabel{4-4}
\| D_t^l w_j(\xi',\cdot,\lambda)\|_{L_2(\R_+)}\le C\left\{
\def\arraystretch{1.2}
\begin{array}{lll}
|\xi'|^{l-j+ \frac 1 2}, & j\le \mu,&l\le \mu,\\
|\xi'|^{1+\mu-j}(\lambda + |\xi'|)^{l-\mu- \frac 1 2},
 & j \le \mu,& l> \mu\\
|\xi'|^{l-\mu}(\lambda + |\xi'|)^{\mu-j+ \frac 1 2},
 & j > \mu,& l\le \mu,\\
(\lambda + |\xi'|)^{l-j+\frac 1 2}, & j > \mu,& l> \mu,
\end{array}\right.
\def\arraystretch{1}
\end{equation}
holds with a constant $C$ not depending on $\xi'$ and $\lambda$.
\end{theorem}

\begin{proof} The existence and the uniqueness of the solution follows
immediately from the fact that $A(\xi',\tau,\lambda)$ (considered as
a polynomial in $\tau$) has exactly $m$ roots with positive imaginary 
part. Let $\gamma(\xi',\lambda)$ be a closed contour in the upper half  of
the complex plane enclosing all roots $\tau_1(\xi',\lambda),\ldots,
\tau_m(\xi',\lambda)$ with positive imaginary part. Then $w_j(\xi',
t,\lambda)$ is given by
\begin{equation}\mylabel{4-5}
w_j(\xi',t,\lambda) = \frac1{2\pi i} \int_{\gamma(\xi',\lambda)}
\frac{M_j(\xi',\tau,\lambda)}{A_+(\xi',\tau,\lambda)} e^{it\tau}\,d\tau
\end{equation}
where 
\begin{equation}\mylabel{4-6}
A_+(\xi',\tau,\lambda) = \prod_{k=1}^m \Big(\tau-\tau_k(\xi',\lambda)\Big)
=: \sum_{k=0}^m a_k(\xi',\lambda)\tau^k 
\end{equation}
and
\begin{equation}\mylabel{4-7}
M_j(\xi',\tau,\lambda) = \sum_{k=0}^{m-j} a_k(\xi',\lambda) \tau^{m-j-k}
\,.
\end{equation}
(Cf., e.g., \cite{agmon-douglis-nirenberg}, Section 1.) The coefficients
are given by the formula of Vieta,
\begin{equation}\mylabel{4-8}
a_k(\xi',\lambda) = \sum_{1\le l_1<\ldots<l_k\le m} (-1)^k \tau_{l_1}
(\xi',\lambda)\cdot\ldots\cdot \tau_{l_k}(\xi',\lambda)\,.
\end{equation}
From \eqnref{4-5} we see, substituting $\tau=r \tilde\tau$, that
\begin{equation}\mylabel{4-9}
 r^{1-j+l}(D_t^l  w_j)(r\xi',\frac t r,r\lambda) = D_t^l w_j(\xi',
t,\lambda)\,,
\end{equation}
and therefore
\begin{equation}\mylabel{4-10}
\| D_t^lw_j(\xi',\cdot,\lambda)\|_{L_2(\R_+)} = r^{\frac 1 2-j+l}
\Big\| D_t^l w_j\Big(\frac{\xi'}{r}, \cdot, \frac\lambda r\Big)
\Big\|_{L_2(\R_+)}\,.
\end{equation}
If we set $r=|\xi'|$ and $\omega' = \frac{\xi'}{|\xi'|}$ we obtain
\begin{equation}\mylabel{4-10a}
\| D_t^lw_j(\xi',\cdot,\lambda)\|_{L_2(\R_+)} = |\xi'|^{\frac 1 2-j+l}
\Big\| D_t^l w_j\Big(\omega', \cdot, \frac\lambda {|\xi'|}\Big)
\Big\|_{L_2(\R_+)}\,.
\end{equation}
The theorem will be proved if we show that for $|\omega'|=1$ we have
\def\arraystretch{1.2}
\begin{equation}\mylabel{4-19}
\|(D_t^l w_j)(\omega',\cdot,\Lambda)\|_{L_2(\R_+)} \le
 \left\{\begin{array}{ll}
 C\,, & j\le \mu\,,\; l\le \mu\,,\\
 C\,\Lambda^{l-\mu-\frac 1 2}, & j\le \mu\,,\; l> \mu\\
 C\,\Lambda^{\mu-j+\frac 1 2}, & j> \mu\,,\; l\le \mu\,,\\
 C\,\Lambda^{l-j+\frac 1 2}, & j> \mu\,,\; l> \mu\,,\\
 \end{array}\right.
\end{equation}
\def\arraystretch{1}
for $\Lambda\ge 1$ and that the left-hand side is bounded by a
constant for $\Lambda\le 1$. 

The boundedness for $\Lambda\le 1$ easily follows from the ellipticity 
of $A(\omega',\tau,\Lambda)$ and inequality \eqnref{1-8}.

For large $\Lambda$ we write 
\[ \gamma(\omega',\Lambda)
=\gamma^{(1)}(\omega',\Lambda)\cup \gamma^{(2)}(\omega',\Lambda)\] 
where $\gamma^{(1)}(\omega',\Lambda)$ encloses 
the zeros $\tau_1(\omega',\Lambda),
\ldots,\tau_\mu(\omega',\Lambda)$ and $\gamma^{(2)}(\omega',\Lambda)$ encloses 
the zeros $\tau_{\mu+1}(\omega',\Lambda)$,
$\ldots,\tau_m(\omega',\Lambda)$.
 Here
we assume that the zeros are numbered according to Lemma \ref{2.3}.
According to this splitting of the contour $\gamma$, we write $w_j
(\omega',t,\Lambda) = w_j^{(1)}(\omega',t,\Lambda) +
w_j^{(2)}(\omega', t, \Lambda)$ with
\[ w_j^{(k)}(\omega',t,\Lambda) := \frac{1}{2\pi i} \int\limits
_{\gamma^{(k)}(
  \omega', t,\Lambda)} \frac{M_j(\omega',\tau,\Lambda)}{ A_+(\omega',
  \tau, \Lambda)} \; e^{it\tau}\, d\tau\quad (k=1,2)\,.\]
 From Lemma \ref{2.3} we know that
\begin{eqnarray*}
|\tau_j(\omega',\Lambda)| & \le & C \qquad\quad(|\omega'|=1,\;
\Lambda\ge\Lambda_0, \quad\quad
j=1,\ldots,\mu)\\
|\tau_j(\omega',\Lambda)| & \le & C\Lambda \qquad\,\,(|\omega'|=1,\; 
\Lambda\ge\Lambda_0,\quad
j=\mu+1,\ldots,m)\,.
\end{eqnarray*}
As $A_{2\mu}$ is elliptic we have, with the notation of Lemma \ref{2.3},
 $|\tau_j(\omega',\Lambda)|\ge C$ for $j=1,\ldots,\mu$ and $|\omega'|=1,\;
\Lambda\ge\Lambda_0$. With our additional assumption we also have
\[ |\tau_j(\omega',\Lambda)|\ge C\Lambda\quad(|\omega'|=1,\; 
\Lambda\ge\Lambda_0,\quad j=\mu+1,\ldots,m),\]
as $\frac{\tau_j(\omega',\Lambda)}{\Lambda}\to \tau_j^1$ and $\Im \tau_j^1
>0$, cf. Lemma \ref{2.3} (ii). Therefore
\def\arraystretch{1.2}
\begin{equation}\mylabel{4-13}
|A_+(\omega',\tau,\Lambda)|= \prod_{k=1}^m |\tau-\tau_k(\omega',\Lambda)|
\ge \left\{\begin{array}{ll} C\Lambda^{m-\mu} & \mbox{ on }\gamma^{(1)}\,,\\
C\Lambda^m& \mbox{ on }\gamma^{(2)}\end{array}\right.
\end{equation}
(note that $|\tau|\approx C$ on $\gamma^{(1)}$ and $|\tau|\approx C\Lambda$
on $\gamma^{(2)}$). Now we have to estimate $|M_j(\omega',\tau,\Lambda)|$ in
\eqnref{4-5}. For this we use the fact that according to \eqnref{4-8}
\begin{equation}\mylabel{4-14}
|a_k(\omega',\Lambda)| \le \sum_{l_1<\ldots<l_k} |\tau_{l_1}|\cdot
\ldots\cdot|\tau_{l_k}|\le \left\{\begin{array}{ll} C\Lambda^k, &
 k\le m-\mu\,,\\
C\Lambda^{m-\mu},& k\ge m-\mu\,.\end{array}\right.
\end{equation}
On $\gamma^{(1)}$ we have 
\[ M_j(\omega',\tau,\Lambda)| \le \left\{\begin{array}{ll} C\Lambda^{m-\mu},
 & j\le \mu\,,\\
C\Lambda^{m-j},& j\ge\mu\,.\end{array}\right.\]
As $\mbox{length}(\gamma^{(1)})\le C$ we obtain
\begin{equation}\mylabel{4-15}
\left| \int_{\gamma^{(1)}} (i\tau)^l \frac{M_j(\omega',\tau,\Lambda)}
{A_+(\omega',
\tau,\Lambda)}e^{i t \tau}\,d\tau\right| \le \left\{\begin{array}{ll}
 C\exp(-Ct)\,, & j\le \mu\,,\\
C\Lambda^{\mu-j}\exp(-Ct),& j\ge\mu\,, \end{array}\right.
\end{equation}
and therefore
\begin{equation}\mylabel{4-16}
\| (D_t^l w_j^{(1)})(\omega',\cdot,\Lambda)\|_{L_2(\R_+)}\le 
 \left\{\begin{array}{ll}
 C\,, & j\le \mu\,,\\
C\Lambda^{\mu-j},& j\ge\mu\,, \end{array}
\right.\quad (|\omega'|=1,\;
\Lambda\ge\Lambda_0)\,.
\end{equation}

For an estimation on $\gamma^{(2)}$ we first remark that for every
$l\ge 0$ we have
\[ |\tau^l M_j(\omega',\tau,\Lambda)| \le \sum_{k=0}^{m-j}
|a_k|\;|\tau^{m-j+l-k}| \le C\Lambda^{m-j+l}\,.\]
Therefore the inequalities
\[ |D_t^l w_j^{(2)}(\omega',t,\Lambda)|\le C\Lambda^{l-j+1}
\exp(-C\Lambda t)\]
and
\begin{equation}\mylabel{4-16a}
\|D_t^l w_j^{(2)}(\omega',\cdot,\Lambda) \|_{L_2(\R_+)}\le C
 \Lambda^{l-j +\frac 1 2} \quad (l\ge 0)
\end{equation}
hold. To find a sharper estimate in the case $j \le \mu$ we use the relation
\begin{eqnarray*}
 \tau^l M_j(\omega',\tau,\Lambda) &=&\tau^{l-j} \sum_{k=0}^{m-j}
a_k(\omega',\Lambda)\tau^{m-k}\\
& =& \tau^{l-j}\Big(A_+(\omega',\tau,\Lambda)-
\sum_{k=m-j+1}^m a_k(\omega',\Lambda)\tau^{m-k}\Big)
\end{eqnarray*}
which yields
\[  D_t^l  w_j^{(2)}(\omega',t,\Lambda) =-\frac{1}{2\pi i}
  \int_{\gamma^{(2)}} \frac{\sum_{k=m-j+1}^m
a_k(\omega',\Lambda)\tau^{m-k+l-j}}{A_+(\omega',
\tau,\Lambda)}e^{i t \tau}\,dt\,.\]
Here we used the fact that the contour $\gamma^{(2)}$ does not enclose the
origin, and therefore $\tau^{l-j}e^{it\tau}$ is holomorphic inside
$\gamma^{(2)}$.
   
We obtain for the case $j\le \mu$ and for every $l\ge 0$ that
\[  \Big|\sum_{k=m-j+1}^m a_k\tau^{m-k+l-j}\Big| \le 
C\Lambda^{m-\mu}\Lambda^{m-
(m-j+1)+l-j}=C\Lambda^{m-\mu+l-1}\]
and
\begin{equation}\mylabel{4-17}
\|D_t^l w_j^{(2)}(\omega',\cdot,\Lambda) \|_{L_2(\R_+)}\le C
 \Lambda^{l-\mu -\frac 1 2} \quad ( j\le\mu,\; l\ge 0)
\end{equation}
in view of Remark \ref{3.0aa} for, say, $\Lambda\ge 1$.
Now we  compare  the right-hand  sides   of
\eqnref{4-16}--\eqnref{4-17}   with the
right-hand side of \eqnref{4-19}.

a) For $j,\, l \le \mu$ the norm of $D^l_tw^{(1)}$ is $O(1)$ and
the norm of $D^l_tw^{(2)}$ is estimated  by $\Lambda^{l-\mu-\frac
1 2}\le\Lambda^{-\frac 1 2}$.
     
b) For $j \le \mu$ and $l>\mu$ according to \eqnref{4-17}
the norm of
$D^l_tw^{(2)}$ is  estimated  by  $\Lambda^{l-\mu-\frac 12}\ge
\Lambda^{\frac 1 2}$ and the norm of $D^l_tw^{(1)}$ is estimated
by a constant.
     
c) For $j > \mu$ and $l\le\mu$ according to \eqnref{4-16} and
\eqnref{4-17}  the norm of
$D^l_tw^{(1)}$ is  estimated  by  $\Lambda^{\mu-j}$ and the norm
of $D^l_tw^{(2)}$ is estimated by $\Lambda^{l-j+\frac 1 2}\le\Lambda^
{\mu-j+\frac 1 2}$.
     
d) For $j,\,l> \mu$ the norm of $D^l_tw^{(2)}$ is  estimated  by
$\Lambda^{l-j+\frac 1 2}$ and the norm of $D^l_tw^{(1)}$ is
estimated by $\Lambda^{\mu-j} <\Lambda^{l-j+\frac 1 2}$.

Thus the inequality \eqnref{4-19} holds, which finishes the 
proof of the theorem.
\def\arraystretch{1}
\end{proof}

\mysection{5}{The main results}

Now we want to prove an a priori estimate for the Dirichlet boundary
value problem corresponding to the elliptic pencil $A(x,D,\lambda)$
defined in \eqnref{1-1}. First we consider model problems in $\R^n$
and $\R^n_+$.

Let $A$ be a polynomial of the form \eqnref{2-1a}. As it was already
mentioned at the beginning of Section 4, the Newton polygon $N_{2m,
  2\mu}$ of $A(\xi,\lambda)$ has the form indicated in Figure 2 with
$r=2m$ and $s=2\mu$. The a priori estimates which we will obtain
below, however, do not use the Sobolev spaces corresponding to this
Newton polygon but the ``energy spaces'' which are defined as the
Sobolev spaces corresponding to the Newton polygon $N_{m,\mu}$. For
this Newton polygon we have
\begin{equation}\mylabel{5-1}
  \Xi(\xi,\lambda) := \Xi_{N_{m,\mu}}(\xi,\lambda) \approx
  (1+|\xi|)^\mu (\lambda+|\xi|)^{m-\mu}\,.
\end{equation}
In the notation of the Introduction, we have $H^\Xi(\R^n) =
H^{(m,\mu)}(\R^n)$. As in Section 2, we will denote by
$\Xi^{(-l)}(\xi,\lambda)$ the weight function corresponding to the
shifted Newton polygon (with a shift of length $l$ to the left). The
space $H^{(-m,-\mu)}(\R^n)$ which appears in the Introduction is equal
to the space $H^{\frac 1\Xi}(\R^n)$.

From the trace results of Theorem \ref{3.3} we immediately obtain the
continuity of the corresponding operators:

\begin{lemma}\mylabel{5.1} {\rm a)} The operator $A(D,\lambda)$ acts
  continuously from $H^\Xi(\R^n)$ to $H^{\frac 1\Xi}(\R^n)$.

\noindent {\rm b)} The boundary operator $D_n^{j-1}$ {\rm (}$j\le
m${\rm )} acts continuously from $H^\Xi(\R^n)$ 
to $H^{\Xi^{(-j+\frac 1 2)}}(\R^{n-1})$.
\end{lemma}

Here and in the following, the continuity of the operator means that
the norm of this operator can be estimated by a constant independent
of $\lambda$.

\begin{proposition}\mylabel{5.2}{\rm (A priori estimate in $\R^n$.)} Let
  $A(\xi,\lambda)$ be $N$-elliptic with parameter in 
  $[0,\infty)$. Then for every 
  $\lambda_0>0$ the inequality  
\begin{equation}\mylabel{5-2} 
  \|u\|_{\Xi,\R^n} \le C\Big( \| A(D,\lambda) u\|_{\frac 1\Xi, \R^n}
+ \lambda^{m-\mu}\| u\|_{L_2(\R^n)}\Big)
\end{equation}
holds for all $\lambda\ge\lambda_0$ with a constant $C=C(\lambda_0)$
independent  of $u$ and $\lambda$.
\end{proposition}

\begin{proof} 
By changing the constant in \eqnref{2-1} we can rewrite the $N$-ellipticity
condition in the form
\begin{eqnarray*}
 \lefteqn{\lambda^{2m-2\mu} +C^{-1}_1\;
 \frac{|A(\xi,\lambda)|^2} {(1+|\xi|^2)^{\mu}
(\lambda^2+|\xi|^2)^{m-\mu}}}\\[2pt]
& &  \ge \lambda^{2m-2\mu}+
|\xi|^{4\mu}\;(1+|\xi|^2)^{-\mu}
(\lambda^2+|\xi|^2)^{m-\mu}\,. 
\end{eqnarray*}
For $|\xi| \ge 1$ the right-hand side can be estimated from below by
\[ (1+|\xi|^2)^{\mu}(\lambda^2+|\xi|^2)^{m-\mu} \]
For $|\xi| \le 1$ and $\lambda \ge \lambda_0$ the right-hand side 
can be estimated from below by
\begin{eqnarray*}
 \lambda^{2m-2\mu}& =& (1+\lambda^{-2})^{-m+\mu}\;
(1+\lambda^{2})^{m-\mu}  \\[2pt]
&  \ge & (1+\lambda_0^{-2})^{-m+\mu}\;
2^{-2\mu}\;(1+|\xi|^2)^{\mu}\;(\lambda^2+|\xi|^2)^{m-\mu} \,.
\end{eqnarray*}
Combining these estimates we obtain for $\lambda \ge \lambda_0$
\[
(1+|\xi|^2)^{\mu}(\lambda^2+|\xi|^2)^{m-\mu} \le C(\lambda_0)
\Big(
 \frac{|A(\xi,\lambda)|^2} {(1+|\xi|^2)^{\mu}{(\lambda^2+
|\xi|^2)^{m-\mu}}}+ \lambda^{2m-2\mu}\Big)\,.\]
Multiplying  both sides by $|F u(\xi)|^2$ and integrating with respect to
$\xi$ we obtain the inequality
\[ \|u\|^2_{\Xi,\R^n} \le C(\lambda_0)\Big( \|A(D,\lambda)u\|^2_{\frac {1}
{\Xi},\R^n} +\lambda^{2m-2\mu}\;\|u\|^2\Big) \]
equivalent to \eqnref{5-2}.
\end{proof}

\begin{theorem}\mylabel{5.3} {\rm (A priori estimate in $\R^n_+$.)} Let
  $A(\xi,\lambda)$ be $N$-elliptic with parameter in $[0,\infty)$ and
  degenerate regularly for $\lambda\to\infty$. Then for every
  $\lambda_0>0$ there exists a constant $C=C(\lambda_0)$ such that for
  all $\lambda\ge\lambda_0$ and all $u\in H^{\Xi}(\R^n_+)$ the
  estimate
\begin{eqnarray}
\|u\|_{\Xi,\R^n_+} &\le & C\Big( \|A(D,\lambda)u\|_{\frac 1
  \Xi,\R^n_+}  \nonumber\\
& + & \sum_{j=1}^m \| D_n^{j-1} u\|_{\Xi^{(-j+\frac 1 2)},\R^{n-1}} +
\lambda^{m-\mu}\|u\|_{L_2(\R^n_+)}\Big)\mylabel{5-7}
\end{eqnarray}
holds.
\end{theorem}

\begin{proof} We will follow a standard plan in elliptic theory. In 
the first part of the proof we reduce \eqnref{5-7}
 to the case $f\equiv 0$.  Then using
Theorem \ref{4.1}, we treat the case of the homogeneous equation.

1) Denote by $E$ a linear operator of extension of functions defined
in $\R^n_+$ to  functions in $\R^n$. If we use the well-known
Hestenes construction then the operator $E: L_2(\R^n_+) \rightarrow
L_2(\R^n)$ and its restriction $E:H^{\Xi}(\R^n_+) \rightarrow     
H^{\Xi}(\R^n)$ are bounded operators. We will denote by $R$ the operator
of restriction of functions on $\R^n$ onto $\R^n_+$.

2) Let $\psi(\xi) \in C^{\infty}(\R^n)$ be a cut-off function, i.e. $\psi
(\xi)=1$ for $|\xi| \le 1$ and $\psi(\xi)=0$ for $|\xi| \ge 2$. We write
\begin{equation}\mylabel{5-8}   
  u = u_1+u_2+v=R\psi(D)Eu+R(1-\psi(D)) A^{-1}(D,\lambda)Ef+v 
\end{equation} 
where we have set $Ef=A(D,\lambda)Eu$. 

First of all we show that $u_1$ and $u_2$ belong to $H^{\Xi}(\R^n_+)$ 
 and their norms in
this space can be estimated by a constant times
\[ \|f\|_{\frac 1 {\Xi},\R^n_+} +\lambda^{m-\mu}\, \|u\|_{L_2(\R^n_+)}\,. \]

3) Since the operator $\psi(D)$ is infinitely smoothing we get for
$\lambda\ge \lambda_0$ that
\[ \|u_1\|_{\Xi,\R^n_+} \le \|\psi(D)Eu\|_{\Xi,\R^n} \le C
\lambda^{m-\mu} \|Eu\|_{L_2(\R^n)}
\le C_1 \lambda^{m-\mu} \|u\|_{L_2(\R^n_+)}\,. \]

4) Using the Fourier transform we obtain
\begin{eqnarray*}
 \|u_2\|_{\Xi,\R^n_+} &\le& 
\|(1-\psi(D))A^{-1}(D,\lambda)Ef\|_{\Xi,\R^n}\\
&=&
\|\Xi(\xi,\lambda)(1-\psi(\xi))A^{-1}(\xi,\lambda)(FEf)
(\xi)\|_{L_2(\R^n)}\,.
\end{eqnarray*}
Since $1-\psi(\xi)=0$ for $|\xi| \le 1$, we obtain from the
$N$-ellipticity  
condition that 
\[ \Xi(\xi,\lambda)\;|1-\psi(\xi)|\;|A^{-1}(\xi,\lambda)| \le C \,
\Xi^{-1}(\xi,\lambda) \]
and
\[\|u_2\|_{\Xi,\R^n_+} \le \mbox{ const }\|Ef\|_{\frac 1 {\Xi},\R^n}\,. \] 
If the norm in $H^{\Xi^{-1}}(\R^n)$ is defined by means of the 
pseudodifferential operator 
\[ \Big( (1+|D'|^2)^{1/2}+iD_n\Big)^{-\mu}\Big((\lambda^2+|D'|^2)^{1/2}
+iD_n\Big)^{-m+\mu}\,, \]
then according to Section 2
\[ \|Ef\|_{\frac 1 {\Xi},\R^n}=\|f\|_{\frac 1 {\Xi},\R^n_+}\,. \]

5) Now we begin the estimation of $v$ defined in \eqnref{5-8}. We have
$v=u-u_1-u_2 \in H^{\Xi}(\R^n_+)$ and
\begin{eqnarray}
\mylabel{5-9} A(D,\lambda)\,v & = & 0 \,,\\
\mylabel{5-10} D_n^{j-1}v(x)|_{x_n=0} & = & h_j(x')\,,
\end{eqnarray}
where we set $h_j(x') := D^{j-1}_nu(x',0)-D^{j-1}_nu_1(x',0)-
D^{j-1}_nu_2(x',0)$.
We shall prove the inequality
\begin{equation}\mylabel{5-11}
\|v\|_{\Xi,\R^n_+} \le \mbox{ const }\Big(\sum_{j=1}^m
\|h_j\|_{\Xi^{(-j+1/2)},\R^{n-1}}
+\lambda^{m-\mu}\|u\|_{L_2(\R^n)}\Big)
\end{equation}
The a priori estimate
\eqnref{5-7} follows from this inequality because, due to Theorem
\ref{3.3}, 
\[ \|D^{j-1}_n u_i\|_{\Xi^{(-j+1/2)},\R^{n-1}} \le \mbox{ const }
\|u_i\|_{\Xi,\R^n_+}  
\quad (i=1,2)\,.\]
The right-hand side of this inequality is already estimated by
the right-hand side of \eqnref{5-7}.

6) We define
\begin{equation}\mylabel{5-12}
\Phi(\xi,\lambda) := \sum_{i,k}|\xi|^i \lambda^k \,,
\end{equation}
where the sum extends over all integer points $(i,k)$ belonging to the
side of $N_{m,\mu}$ which is not parallel to the coordinate lines.
From this definition it follows that
\begin{equation}\mylabel{5-13}
\Phi(\xi,\lambda)\approx |\xi|^{\mu} (\lambda+|\xi|)^{m-\mu}\,.
\end{equation}
and $\|v\|_{\Xi,\R^n_+}$ is equivalent to
\[ \|v\|_{\Phi,\R^n_+}+ \lambda^{m-\mu} \|v\|_{L_2(\R^n_+)}\,.\]
The second term can be estimated by $\lambda^{m-\mu}(\|u\|_{L_2(\R^n_+)} +
\|u_1\|_{L_2(\R^n_+)}  + \|u_2\|_{L_2(\R^n_+)}) \le 
\lambda^{m-\mu}\|u\|_{L_2(\R^n_+)}+ \|u_1\|_{\Xi,\R^n_+}+
\|u_2\|_{\Xi,\R^n_+}$. Therefore,
it is enough to estimate $\|v\|_{\Phi,\R^n_+}$ 
by the right-hand side of \eqnref{5-11}.
Repeating the argument in Section 2 (see \eqnref{3-10a}) 
we reduce our problem
to the estimation of
\[ \int_0^{\infty}\|(D^l_nv)(\cdot,x_n)\|^2_{\Phi^{(-l)},\R^{n-1}}\;
dx_n
\quad (l=0,\ldots,m) \]
or after the Fourier transform with respect to $x'$
\begin{equation}\mylabel{5-14}
 \int_0^{\infty}\int_{\R^{n-1}} |\Phi^{(-l)}(\xi,\lambda)(D^l_nF'v)
(\xi',x_n)|^2\, d\xi' d x_n \quad (l=0,\dots,m) \,.
\end{equation}
The function $F'v(\xi',x_n)=:w(\xi',x_n)$ is (for almost every 
$\xi'\in\R^{n-1}$) a solution of
\begin{eqnarray}
\mylabel{5-15} A(\xi',D_n,\lambda) w(x_n) & = & 0 \,,\\
\mylabel{5-16} D_n^{j-1}w(x_n)|_{x_n=0} & = & (F'h_j)(\xi')\,.
\end{eqnarray}
Due to Theorem \ref{4.1}, this solution is unique and given by
\begin{equation}\mylabel{5-17}
w(\xi',x_n) = \sum_{j=1}^m w_j(\xi',x_n,\lambda) (F'h_j)(\xi')
\end{equation} 
with $w_j(\xi',x_n,\lambda)$ being the solution of
\eqnref{4-2}--\eqnref{4-3}. 

8) To obtain the estimate for $w = F'v$ we reformulate Theorem \ref{4.1}. 
It follows from the definition of $N_{m,\mu}$ that 
\def\arraystretch{1.2}
\[ \Phi^{(-r)}(\xi,\lambda) \le \left\{ \begin{array}{ll}
 |\xi|^{\mu-r}(\lambda+|\xi|)^{m-\mu},& r \le \mu,\\
 (\lambda+|\xi|)^{m-r},& r>\mu.\end{array}\right. \]
From this it follows that
\[ \frac {\Phi^{(-j+1/2)}(\xi',\lambda)}{\Phi^{(-l)}(\xi',\lambda)} \le
\left\{
\begin{array}{lll}
 C|\xi'|^{l-j+\frac 1 2}, &  l\le\mu, &j\le\mu,\\
 C|\xi'|^{\mu-j+\frac 1 2}(\lambda+|\xi'|)^{l-\mu},& l>\mu,& j\le\mu,\\
 C|\xi'|^{l-\mu}(\lambda+|\xi'|)^{\mu-j+\frac 1 2},& l \le \mu, &j>\mu,\\
 C(\lambda+|\xi'|)^{l-j+\frac 1 2},& l > \mu, &j>\mu.
\end{array}\right.\]
Comparing the right-hand sides of these inequalities with the 
right-hand side of \eqnref{4-4} we see that
\[ \| D_n^l w_j(\xi',x_n,\lambda)\|_{L_2(\R_+)} \le C
\frac {\Phi^{(-j+1/2)}(\xi',\lambda)}{\Phi^{(-l)}(\xi',\lambda)}\,. \]
From \eqnref{5-17} and the last inequality it follows that
\begin{eqnarray*}
\lefteqn{\hspace*{-1.5cm} (\Phi^{(-l)}(\xi,\lambda))^2\int_0^{\infty}
|D_n^l w(\xi',x_n,\lambda)|^2\,dx_n }\\
&\le &
C\; \sum|\Xi^{(-j+\frac 1 2)}(\xi',\lambda)(F'h_j)(\xi')|^2 \,.
\end{eqnarray*}
Integrating this inequality with respect to $\xi'$ we obtain 
the desired estimate.
\end{proof}

Now we consider the Dirichlet boundary value problem for
differential operators with parameter acting
on a smooth compact manifold $M$ with  smooth boundary $\Gamma$. In
this case we can choose a finite number of coordinate systems. In each
of these systems the operator is of the form \eqnref{1-1}. The principal part 
of the operator is invariantly defined at each of these systems and at 
every fixed point $x^0 \in M$ it is of the form
\begin{equation}\mylabel{5-29}
A^{(0)}(x^0,D,\lambda) = A^{(0)}_{2m}(x^0,D) + \ldots + \lambda^{2m-2\mu}
A^{(0)}_{2\mu}(x^0,D)
\end{equation}
(here $A^{(0)}_j$ denotes the principal part of $A_j$). We suppose that for 
each $x^0 \in \overline M$ our operator is $N$-elliptic with parameter.
From the reason of continuity and compactness the constant $C$ in inequality
\eqnref{2-1} can be chosen independent of $x^0$.

We can suppose without loss of generality that the coefficients 
of $A(x,D,\lambda)$ are of the form
\begin{equation}\mylabel{eq5.15}
a_{\alpha j}(x)= a_{\alpha j}+a'_{\alpha j}(x),\quad a_{\alpha j}
\in {\cal D}(\R^n)  \,.
\end{equation}

Now we fix a point $x^0 \in \Gamma$ and a coordinate system in the 
neighborhood of $x^0$ such that in this system locally the boundary $\Gamma$
is given by the equation $x_n=0$. In this case we can define an analog
of the polynomial \eqnref{2-7}:
\begin{equation}\mylabel{5-29a}
Q(x^0,\tau)=\tau^{-2\mu}A^{(0)}(x^0,0,\tau,1)     
\end{equation}
Suppose that at a point $x^0 \in \Gamma$ and in a fixed coordinate system
this polynomial has $m-\mu$ roots in the upper half-plane of the complex 
plane. It easily follows from this fact that every polynomial
\eqnref{5-29a}
 corresponding
to an arbitrary $x^0\in\Gamma$ has the same property.
 In this case we say that the
operator $A(x,D,\lambda)$ degenerates regularly at the boundary $\Gamma$.

\begin{lemma}\mylabel{5.5} For $a(x) = a+a'(x)$ with $a'\in 
{\cal D}(\R^n)$ and $f\in
    H^{\frac 1 \Xi}(R^n)$ we have $af\in H^{\frac 1 \Xi}(R^n)$, and
the following statements hold:

\noindent {\rm a)} There exists a constant $C(a)$ depending on $a$ but 
not on $f$ or $\lambda$ such that
\begin{equation}\mylabel{5-30}
\| af\|_{\frac 1\Xi,\R^n} \le C(a) \|u\|_{\frac 1
  \Xi,\R^n}\,. 
\end{equation}

\noindent {\rm b)} There exists a constant $C'(a)$ depending only on
$a$ such that the inequality
\begin{equation}\mylabel{5-31}
  \| af\|_{\frac 1\Xi,\R^n} \le \sup_{x\in\R^n} |a(x)| \|f\|_{\frac 1
    \Xi, \R^n} + C'(a) \|f\|_{\Psi,\R^n}
\end{equation}
holds, where we have set
\begin{equation}\mylabel{5-32}
\|f\|_{\Psi,\R^n} := \left(\int (1+|\xi|)^{-2\mu-2}
  (\lambda+|\xi|)^{-2m +2\mu}|\hat f(\xi)|^2 \, d\xi\right)^{\frac 1
  2}\,.
\end{equation}
\end{lemma}

\begin{proof} Part a)  is a special case of the following 
more general result which is
taken from \cite{volevich-panejah}, Section I.2.4. Let $\sigma$ be a weight
function which satisfies 
\[ \sigma(\xi)\sigma^{-1}(\eta) \le C(1+|\xi-\eta|^m)\,.\]
Then we have for $a'\in {\cal D}(\R^n)$ the inequality
\begin{equation}\mylabel{5-35}
 \| a'f\|_{H^\sigma(\R^n)} \le c(a') \|f\|_{H^\sigma(\R^n)}
\end{equation}
with $c(a') := C \int (1+|\xi|^m) |(F a')(\xi)| d\xi$. 

Part b) can be shown by standard arguments similar to those used in
 \cite{lions-magenes}, Section 1.7.1, and \cite{gindikin-volevich},
 Lemma 1.4.5.
\end{proof}

Using the above mentioned covering of $M$ by local coordinate systems and a 
partion of unity subordinated to this covering we can define the spaces
$H^{\Xi}, H^{\frac {1}{\Xi}}$ and $H^{\Xi^{(-j+3/2)}}$. From Lemma
\ref{5.5} and the
trace results for model problems in $\R^n$ and $\R^n_+$ we immediatly obtain

\begin{lemma}\mylabel{5.6} Let $D_\Gamma(u) := (u|_\Gamma,
\frac\partial{\partial\nu}u|_\Gamma,\ldots,(\frac\partial{\partial\nu}
)^{m-1} u|_\Gamma)$ be the Dirichlet boundary operator. Then
\[ (A(x,D,\lambda),D_\Gamma):\; H^{\Xi}(M) \longrightarrow
 H^{\frac 1\Xi}(M) \times
\prod_{j=1}^m H^{\Xi^{(-j+\frac 1 2)}}(\Gamma)\]
is continuous with norm bounded by a constant independent of
$\lambda$.
\end{lemma}
  
\begin{theorem}\mylabel{5.7} Let $A(x,D,\lambda)$ be an operator
  pencil of the form \eqnref{1-1}, acting on the manifold $M$ with
  boundary $\Gamma$. Let $A$ be $N$-elliptic with parameter in
  $[0,\infty)$ and assume that $A$ degenerates regularly at the
  boundary $\Gamma$. Then for $\lambda\ge\lambda_0$ there exists a
  constant $C=C(\lambda_0)$, independent of $u$ and $\lambda$, such
  that
\begin{eqnarray}
\mylabel{5-40}
\|u\|_{\Xi,M} &\le& C\Big( \|A(x,D,\lambda)u\|_{\frac 1\Xi,M} +
\sum_{j=1}^m \Big\| \Big(\frac{\partial}{\partial\nu}\Big)^{j-1}u
\Big\|_{\Xi^{(-j+\frac 1 2)},\Gamma}\\
& & +\quad \lambda^{m-\mu} \|u\|_{L_2(M)}\Big)\,.\nonumber
\end{eqnarray}
\end{theorem}

\begin{proof} For the proof we use the standard technique of
  localization (``freezing the coefficients''). We only indicate the
  main steps. By means of a partition of unity it is sufficient to
  prove \eqnref{5-40} for $u\in H^\Xi(M)$ with small support ${\rm
    supp\, } u\subset U$. In the case $U\cap \Gamma =\emptyset$, we fix
  $x_0\in U$ and use local coordinates. We obtain from the a priori
  estimate for the model problem in $\R^n$ that
\begin{eqnarray}
\|u\|_{\Xi,\R^n} & \le & C_1\Big( \| A^{(0)}(x_0,D) u\|_{\frac 1
  \Xi,\R^n} + \lambda^{2m-2\mu} \|u\|_{L_2(\R^n)}\Big) \nonumber\\
& \le & C_1\Big( \| A(x,D) u\|_{\frac 1
  \Xi,\R^n} + \lambda^{2m-2\mu} \|u\|_{L_2(\R^n)}\Big) \nonumber\\
& & +\, C_1 \| (A(x,D) - A^{(0)}(x_0,D))u\|_{\frac 1
  \Xi,\R^n} \mylabel{5-41}
\end{eqnarray}
with a constant $C_1$ independent of $u$ and $\lambda$.

We fix $\epsilon > 0$. From Lemma \ref{5.5} b) we obtain if the
support of $u$ is sufficiently small that 
\begin{equation}\mylabel{5-42}
\| (A(x,D) - A^{(0)}(x_0,D))u\|_{\frac 1
  \Xi,\R^n} \le \epsilon \|u\|_{\Xi,\R^n} + C\|u\|_{\Xi^{(-1)},\R^n}\,.
\end{equation}
Here we have used that
\begin{equation}\mylabel{5-43}
 \sum_{\alpha,k} \lambda^k \| D^\alpha u\|_{\frac 1 \Xi,\R^n} \le C
 \|u\|_{\Xi,\R^n}
\end{equation}
and
\begin{equation}\mylabel{5-44}
\sum_{\alpha,k} \lambda^k \| D^\alpha u\|_{\Psi,\R^n} \le C
 \|u\|_{\Xi^{(-1)},\R^n}
\end{equation}
Now we use the interpolation inequality
\begin{equation}\mylabel{5-45}
\|u\|_{\Xi^{(-1)},\R^n}\le \epsilon \|u\|_{\Xi,\R^n} +
C\lambda^{m-\mu} \|u\|_{L_2(\R^n)}
\end{equation}
which is a consequence of the interpolation inequality for the Sobolev 
spaces $H^s(\R^n)$ because of
\begin{equation}\mylabel{5-46}
\|u\|_{\Xi^{(-1)},\R^n}\approx \|u\|_{H^{m-1}(\R^n)} + \lambda^{m-\mu} 
\|u\|_{H^{\mu-1}(\R^n)}\,.
\end{equation}
If we choose $\epsilon$ with $C_1\epsilon<1$ we obtain
\begin{equation}\mylabel{5-47}
\|u\|_{\Xi,\R^n} \le C\Big( \|A(x,D,\lambda)u\|_{\frac 1\Xi,\R^n} +
\lambda^{m-\mu} \|u\|_{L_2(\R^n)}\Big)\,.
\end{equation}
In the case $U\cap \Gamma\not=\emptyset$ we choose $x_0\in U\cap
\Gamma$, use local coordinates, and obtain in the same way as above
\begin{eqnarray}
\|u\|_{\Xi,\R^n_+} &\le&  C\Big( \|A(x,D,\lambda)
u\|_{\frac 1\Xi,\R^n_+}
+ \sum_{j=1}^m \| D_n^{j-1} u\|_{\Xi^{(-j+\frac 1 2)},\R^{n-1}} \nonumber\\
& & \qquad +\; \lambda^{m-\mu} \|u\|_{L_2(\R^n)}\Big)\,,\mylabel{5-48}
\end{eqnarray}
where we used the a priori estimate for $(A^{(0)}(x_0,D),
(D_n^{j-1})_{j=1}^m)$.
\end{proof}

Now to finish the existence theory we present the construction of the
right (rough) parametrix of the Dirichlet problem.

Suppose the assumptions of Theorem \ref{5.7} to hold. We will see
below that  the
solution constructed above for constant coefficients is a right
parametrix with respect to the Sobolev spaces defined by the
Newton polygon. We will 
define this parametrix, as usual, with the help of local coordinates.
 First  of  all  we  present  the
construction in the case of model domains $\R^n$ and $\R^n_+$.

\begin{proposition}\mylabel{5.7a} Suppose 
$A(x,\xi,\lambda)$ satisfies the $N$-ellipticity condition and the
coefficients are of the form \eqnref{eq5.15}.
Then there exists a bounded operator
\begin{equation}\mylabel{eq29}
B: H^{\frac {1}{\Xi}}(\R^n) \rightarrow H^{\Xi}(\R^n) 
\end{equation}
such that
\begin{equation}\mylabel{eq30}
A(x,D,\lambda)B=I+T     
\end{equation}
where $I$ denotes the identity operator in $H^{\frac1{\Xi}}
(\R^n)$ and
\begin{equation}\mylabel{eq31}
T: H^{\frac 1{\Xi}}(\R^n) \to H^{\Theta}(\R^n)   
\end{equation}
is continuous with norm bounded by a constant independent of
$\lambda$. Here we posed $\Theta(\xi,\lambda)=(1+|\xi|)/
\Xi(\xi,\lambda)$.
\end{proposition}

\begin{proof} We define $B$ as a classical ps.d.o with symbol
\begin{equation}\mylabel{eq32}
 B(x,\xi,\lambda):= \psi(\xi) \frac 1{A_0(x,\xi,\lambda)},
\end{equation}                                    
where $\psi\in C^\infty(\R^n)$ is a cut-off function with $\psi\equiv 0$
for $|\xi|\le 1$ and $\psi\equiv 1$ for $|\xi|\ge 2$.

The continuity of operator \eqnref{eq29} is equivalent to the statement that
the $L_2-L_2$ norm of the operator
\[ (1+|D|^2)^{\frac{\mu}{2}}
(\lambda^2+|D|^2)^{\frac{m-\mu}{2}} B(x,D,\lambda) (1+|D|^2)^{\frac
{\mu}{2}}(\lambda^2+|D|^2)^{\frac{m-\mu}{2}} \]
can be estimated by a constant independent of $\lambda$. Using
standard results on the $L_2$-boundedness of
ps.d.o. (cf. \cite{kumano-go}, Section 2.4) we have to show the
inequalities
\[ \psi(\xi)\Big|D^{\alpha}_x A_0^{-1}(x,\xi,\lambda)\Big| \le C_{\alpha}
(1+|\xi|)^{-2\mu}(\lambda+|\xi|)^{-2m+2\mu}\,. \]
For $|\alpha|=0$ this inequality directly follows from $N$-ellipticity 
with parameter, to prove it for arbitrary $\alpha$ we must
use the chain rule.

To prove \eqnref{eq30}--\eqnref{eq31} we write the operator $T$ in the form
\[ T = \tilde T + (A(x,D,\lambda)-A_0(x,D,\lambda)) B\]
with $\tilde Tu = A_0(x,D,\lambda) B u-u$. Noting that
\[ A(x,D,\lambda)-A_0(x,D,\lambda) : H^\Xi(\R^n)\to H^{\Theta}(\R^n)\]
is continuous, it is sufficient to prove \eqnref{eq31} with $T$
replaced by  $\tilde T$.  As above,  this is equivalent to the
uniform $L_2 - L_2$ boundedness of
\[ (1+|D|^2)^{-\frac{\mu-1}{2}}
(\lambda^2+|D|^2)^{-\frac{m-\mu}{2}} \tilde T \;  (1+|D|^2)^{\frac{\mu}{2}}
(\lambda^2+|D|^2)^{\frac{m-\mu}{2}}\,. \]
For this it is enough to show that the symbol $\tilde T(x,\xi,\lambda)$
of $\tilde T$  satisfies
\begin{equation}\mylabel{5-56}
(1+|\xi|) \Big| D_x^\beta \tilde T(x,\xi,\lambda)\Big| \le C_\beta \,.
\end{equation}
The last inequality follows easily from the fact that for $|\xi|\ge 2$
we have
\begin{equation}\mylabel{5-57}
\tilde T(x,\xi,\lambda) = \sum_{0<|\alpha|\le 2m} \frac{1}{\alpha!}
 \partial_\xi^\alpha A_0(x,\xi,\lambda)\;D_x^\alpha
 \;\frac{1}{A_0(x,\xi,\lambda)}
\end{equation}
and from the estimates
\begin{equation}\mylabel{5-58}
\Big| D_x^\beta \; \frac 1{A_0(x,\xi,\lambda)} \Big| \le C\;
(\Xi_P(\xi,\lambda))^{-1} \quad (|\xi|\ge 2)
\end{equation}
and
\begin{equation}\mylabel{5-59}
|D^{\gamma}_x\partial_\xi^\alpha A_0(x,\xi,\lambda)| \le C \;
\Xi_P^{(-|\alpha|)}(\xi,\lambda) \quad (0\le |\alpha|\le 2m)\,.
\end{equation}
\end{proof}

\begin{proposition}\mylabel{5.8}
Suppose the conditions of Theorem {\rm \ref{5.7}} are satisfied
and the coefficients of $A(x,\xi,\lambda)$  are  of  the  form
\eqnref{eq5.15}. Then there exists a bounded operator
\begin{equation}\mylabel{eq5.41}
B: H^{\frac 1{\Xi}}(\R^n_+) \times \prod_{j=1}^m H^{\Xi^{(-j+1/2)}}
(\R^{n-1}) \to H^{\Xi}(\R^n_+)
\end{equation}               
such that
\begin{equation}\mylabel{5.38a}
(A,\gamma_0, \gamma_0 D_n, \dots,\gamma_0 D_n^{m-1})^t B = I + T \,,
\end{equation}
where $I$ denotes the identity operator in 
\[ H^{\frac1{\Xi}}(\R^n_+)
\times\prod_{j=1}^m H^{\Xi^{(-j+1/2)}}(\R^{n-1})\,,\]
 $\gamma_0$ is the
operator of taking the trace of the function at ${x_n=0}$, and
\begin{equation}\mylabel{eq5.42}
T: H^{\frac 1{\Xi}}(\R^n_+) \times \prod_{j=1}^m H^{\Xi^{(-j+1/2)}}
(\R^{n-1}) \to H^{\Theta}(\R^n_+) \times \prod_{j=1}^m H^{\Xi^
{(-j+3/2)}}(\R^{n-1})
\end{equation}                  
is continuous with norm bounded by a constant independent of
$\lambda$.
\end{proposition}

\begin{proof} We set
\begin{equation}\mylabel{eq5.43}
B(f,g_1,\ldots,g_m):=B_0 f+ \sum_{j=1}^m B_j(g_j-\gamma_0
 D_n^{j-1} B_0f)\,.
\end{equation}                         
Here (compare Proposition \ref{5.7a})
\begin{equation}\mylabel{5-50}
 B_0 f := R \psi(D) A_0^{-1}(x,D,\lambda)Ef
\end{equation}
for $f\in H^{\frac 1\Xi}(\R^n_+)$. As in \eqnref{5-8} $R$ and $E$ are
restriction and extension operators, respectively,  $\psi(\xi)$
is the cut-off function from the proof of
Proposition \ref{5.7a} and $B_j$ for $j=1,
\dots,m$ is a ps.d.o. in $\R^{n-1}$ (with $x_n$ as parameter) 
defined by
\begin{equation}\mylabel{5-51}
(B_j g_j)(x',x_n) :=\psi'(D')w_j(x',x_n,D',\lambda) g_j \,.
\end{equation}
The symbol $w_j$ in \eqnref{5-51} is given by (compare \eqnref{4-5})
\begin{equation}\mylabel{5-52}
 w_j(x',x_n,\xi',\lambda) := \frac{1}{2\pi i}
 \int_{\gamma(\xi',\lambda)} \frac{M_j(x',0,\xi',\tau,\lambda)}{
   A_+(x',0,\xi',\tau,\lambda)} e^{ix_n\tau}\,d\tau\,,
\end{equation}
where $A_+(x',0,\xi',\tau,\lambda)$ and $M_j(x',0,\xi',\tau,\lambda)$
are given by \eqnref{4-6}--\eqnref{4-7} with $A(\xi',\lambda)$
replaced by $A^{(0)}(x',0,\xi',\lambda)$. The function
$\psi'(\xi') \in C^\infty(\R^{n-1})$ is defined by $\psi'(\xi') :
= \psi(\xi',0)$.

First of all we check that the operator \eqnref{eq5.43} is bounded. The
boundness of   $B_0$   follows   easily   from  the  proof  of
Proposition \ref{5.7a}. Now we proof the continuity of operators
\begin{equation}\mylabel{5-60}
B_j: H^{\Xi^{(-j+1/2)}}(\R^{n-1}) \to H^\Xi(\R^n_+)\,.
\end{equation}
If the condition \eqnref{eq5.15} is satisfied, then the operator $B_j$
can be represented in the form $B_j =B^0_j+B'_j$, where $B^0_j$
is a ps.d.o with symbol independent of $x$ and the symbol of
$B'_j$ is independent of $x$, when the modulus of $x$ is large enough.

In fact, using the norm
\begin{equation}\mylabel{5-61}
\Big( \sum_{l=0}^m \int_0^\infty \| D_n^l
u(\cdot,x_n)\|^2_{\Xi^{(-l)}, \R^{n-1}}\, d x_n\Big)^{1/2}
\end{equation}
in $H^\Xi(\R^n_+)$ and Theorem \ref{4.1} we can  easily  prove  the
continuity of operator $B^0_j$. To prove the continuity of
$B'_j$ we show the  estimates
\begin{equation}\mylabel{5-61a}
\psi'(\xi')\Big(\int_0^{\infty} \Big| D_{x'}^{\beta'}D^l_n w_j(x',
x_n,\xi',\lambda)\Big|^2\,dx_n\Big)^{1/2} \le C_{\beta'} \; \frac{\Xi^{(-j+1/2)}
(\xi',\lambda)}{\Xi^{(-l)}(\xi',\lambda)}\,.
\end{equation}
The case $\beta'=0$ was treated in the proof of Theorem \ref{4.1},
the general case follows by the same method after
differentiating in \eqnref{5-52} under the integral sign.

Now we prove the continuity of the operator \eqnref{eq5.42}. The main step of
the proof is reduced to showing that $AB_j=0,\quad j=1,\dots,m,$
up to a lower order operator, i.e. the operator
\begin{equation}\mylabel{eq5-45}
A(x,D,\lambda)B_j: H^{\Xi^{(-j+1/2)}}(\R^{n-1}) \to H^\Theta(\R^n_+)
\end{equation}
is continuous.

Before proving  this  statement we finish the proof of \eqnref{5.38a}.
Denote by $T_0$, $T_1,\dots, T_m$ the components of operator $T$.
Operator $T_0$ maps the space on the
 left-hand side of  \eqnref{eq5.42} into
$H^{\Theta}(\R^n_+)$ and $T_k,\;  k>0,$ maps the space on the
left-hand  side of
\eqnref{eq5.42} into $H^{\Xi^{(-k+3/2)}}(\R^{n-1})$. We have
\[ T_0\Big\{f,g_1, \dots, g_m\Big\}= (AB_0f-f)+\sum_{j=1}^m AB_j
(g_j-\gamma_0 D_n^{j-1}B_0f)\,.\]
The first   term   on   the   right-hand   side   belongs   to
$H^{\Theta}(\R^n_+)$ according to  Proposition  \ref{5.7a},  the  sum
belongs to this space according to \eqnref{eq5-45}.

Turning to estimates of other components of $T$ we note that
according to \eqnref{4-3}
\[ \gamma_0 D_n^{k-1}B_j(x',x_n,D_n,\lambda)h=\delta_{jk}\psi'(D')h\,. \]
From this follows that for $k \ge 1$ we have
\begin{eqnarray}
\lefteqn{T_k\Big\{f,g_1,\dots,g_m\Big\}}\nonumber\\
&=&\gamma_0 D^{k-1}_n B_0f+\sum_{j=1}^m
\gamma_0D^{k-1}_n B_j(g_j-\gamma_0 D^{k-1}_n B_0 f)-g_k\nonumber\\
&=&(1-\psi'(D'))(\gamma_0 D^{k-1}_n B_0f-g_k) \,.\mylabel{eq5.47}
\end{eqnarray}
The function $1-\psi'(\xi')$ belongs to ${\cal D}(\R^{n-1})$
and, consequently, function \eqnref{eq5.47}
 belongs to $H^{\infty}(\R^{n-1})$ 
for arbitrary $\{f,g_1,\dots,g_m\}$ from the  space on the left-hand
side of \eqnref{eq5.42}.

To prove \eqnref{eq5-45} we can suppose, without loss of generality, that
the operator $A$ is replaced by its principal part
\[ A^{(0)}(x,D,\lambda)=A^{(0)}(x',0,D,\lambda)+\Big(A^{(0)}(x,D,\lambda)-
A^{(0)}(x',0,D,\lambda)\Big)\,. \]
The composition of $A^{(0)}(x',0,D,\lambda)$ and $B_j$ is a ps.d.o. in
$\R^{n-1}$ (depending on $x_n$) with the symbol 
\begin{eqnarray*}
\lefteqn{C_j(x',x_n,\xi',\lambda)}\\
&=&  \psi'(\xi')\sum_{|\alpha'|=1}^{2m}
\frac{1}{(\alpha')!} \partial_{\xi'}^{\alpha'}
A^{(0)}(x',0,\xi',D_n,\lambda) D_{x'}^{\alpha'} w_j(x',x_n,\xi',\lambda)\,.
\end{eqnarray*}
The summand corresponding to $\alpha'=0$ is identically zero due
to the definition of  $w_j$. We can rewrite the symbol above in the form
\[
\sum_{l=0}^{2m-1} c_l(x',\xi',\lambda)D_n^l w_j(x,\xi',\lambda),
\]
where the coefficients $c_l(x',\xi',\lambda)$ and their derivatives with 
respect to $x'$ admit the estimates
\[ 
\Big|D^{\beta'}_{x'}c_l(x',\xi',\lambda)\Big|\le C\;
\Xi_P^{(-l-1)}(\xi',\lambda)\,.
\]
Differentiating under the integral sign and repeating the argument of
Theorem \ref{4.1} we obtain
\[
\Big(\int_0^{\infty}\Big |D^{\beta'}_{x'}D_n^lw_j(x',x_n,\xi',\lambda)
dx_n \Big|\Big)^{1/2} \le \frac{\Xi^{(-j+\frac {1}{2})}(\xi',\lambda)}
{\Xi^{(-l)}(\xi',\lambda)}\,.
\]
Remembering that $\Xi_P(\xi,\lambda) \equiv \Xi^2(\xi,\lambda)$,
 we easily deduce
\[
\Big(\int_0^{\infty}\Big |D^{\beta'}_{x'}C(x',x_n,\xi',\lambda)
\Big|^2dx_n \Big)^{\frac {1}{2}} \le C\;(1+|\xi'|)^{-1}\Xi^{(-j-
\frac{1}{2})}(\xi',\lambda)\;\Xi(\xi',\lambda)\,.
\]
From this follows the continuity of operator \eqnref{eq5-45} with $A$ replaced 
by $A^{(0)}(x',0,D,\lambda)$.

Applying the Lagrange formula to the coefficients of the operator
$A^{(0)}(x,D,\lambda)-A^{(0)}(x',0,D,\lambda)$, we can rewrite this operator 
in the form 
\[
\sum_{l=0}^{2m}c'_l(x,\xi',\lambda)x_n D^l_n\,,
\]
where
\[
|D^{\beta'}_{x'} c'_l(x,\xi',\lambda)| \le C\;\Xi(\xi',\lambda)\;
\Xi^{(-l)}(\xi',\lambda)\,.
\]
Using the relation
\[
x_n\exp(ix_n\tau)=-i\frac{\partial}{\partial \tau}\exp(ix_n\tau)
\]
and integrating under the sign of the contour integral in \eqnref{4-5}
 we obtain that
\[
x_n D^l_n w_j(x',x_n,\xi',\lambda)=il\; D^{l-1}_n w_j(x',x_n,\xi',\lambda).
\]
Now we easily come to the inequalities
\[
\Big (\int_0^{\infty}\Big |x_n D^{\beta'}_{x'}D_n^lw_j(x',x_n,\xi',
\lambda)\Big|^2 dx_n \Big)^{\frac {1}{2}} \le C\;\frac
{\Xi^{(-j+\frac {1}{2})}(\xi',\lambda)}{(1+|\xi'|)
 \Xi^{(-l)}(\xi',\lambda)}\,.
\]
The above estimates permit us to finish the proof of the Proposition.
\end{proof}

Now we return to the case of a manifold as considered in Theorem
\ref{5.7}. To construct a right parametrix for $(A,D_\Gamma)$, we will 
use local coordinates.

 Let $\{{\cal O}_j\}_{j=1,\ldots,N}$ 
be a covering of $M$ with
coordinate neighborhoods where ${\cal O}_j$ is homeomorphic to an open 
ball in $\R^n$ for $j=1,\ldots,N'$ and homeomorphic to a semi-ball
$\{x\in\R^n: |x|<r_j,\; x_n>0\}$ for $j=N'+1,\ldots,N$. In local
coordinates corresponding to the neighborhood ${\cal O}_k$ we obtain
an operator $A_k$ in $\R^n$ for $k=1,\ldots,N'$ and in $\R^n_+$ for
$k=N'+1,\ldots, N$. We may assume that the coefficients of $A$ in
these local coordinates have the form \eqnref{eq5.15}.

For $k=1,\ldots,N'$ we define a local right parametrix $B_k$ to the operator 
$A_k$ according to Proposition \ref{5.7a}; for $k=N'+1,\ldots,N$ this
can be done due to Proposition \ref{5.8}. We then set
 \[ B :=\sum_{j=1}^N \psi_k B_k(\phi_k\; \cdot)\]
 where $\{\phi_k\}_k$ is a
partition of unity subordinated to the covering of $M$ and $\psi_k$
has support in ${\cal O}_k$ and is equal to 1  in a neighborhood of
${\rm supp\, } \phi_k$. From the previous two propositions we obtain
the following theorem.

\begin{theorem}\mylabel{5.9} The operator $B$ defined above is a right rough
  parametrix to $(A,D_\Gamma)$ in the sense that
\begin{equation}\mylabel{5-53}
(A,D_\Gamma) B = I + T \,,
\end{equation}
where $I$ denotes the identity operator in $H^{\frac 1{\Xi}}(M) \times
\prod_{j=1}^m H^{\Xi^{(-j+1/2)}}(\Gamma)$ and 
\begin{equation}\mylabel{5-54}
T: H^{\frac 1{\Xi}}(M) \times \prod_{j=1}^m H^{\Xi^{(-j+1/2)}}(\Gamma) \to
H^{\Theta}(M) \times \prod_{j=1}^m H^{\Xi^{(-j+3/2)}}(\Gamma)
\end{equation}
is continuous with norm bounded by a constant independent of
$\lambda$.
\end{theorem}

As the norms considered in \eqnref{5-54} are for every fixed
$\lambda$ equivalent to the corresponding norms in the standard
Sobolev spaces, we obtain the
compactness of $T$ and the Fredholm property of $(A,D_\Gamma)$. Note
also that for fixed $\lambda$ the operator $(A,D_\Gamma)$ is elliptic
in the usual sense. From the last remark we also see that the index of 
this operator is equal to zero.

\newpage

\def\refname

\end{document}